\documentclass{article}
\usepackage[T1,T2A]{fontenc} %
\usepackage[cp1251]{inputenc} %
\usepackage[russian]{babel} %
\usepackage{amsfonts} %
\usepackage{amsmath, bm} %
\usepackage{amssymb} %
\usepackage{MnSymbol}
\usepackage{color}
\usepackage[active]{srcltx}
\usepackage{xr}
\usepackage[final]{pdfpages}
\usepackage{graphicx}
\usepackage{wrapfig}
\usepackage{tikz}
\usetikzlibrary{arrows,decorations.pathmorphing,
backgrounds,positioning,fit,petri}
\usepgflibrary{arrows} 
\usetikzlibrary{arrows} 
\definecolor{whitesmoke}{rgb}{0.9,0.9,0.9}
\definecolor{whitesmoke2}{rgb}{0.8,0.8,0.8}
\usetikzlibrary{backgrounds}
\begin{document}

\noindent 
╙─╩ 517.911 \\ 
MCS 34A12 

\smallskip
{\large
\noindent 
{\bf ╬ ёє∙хёЄтютрэшш Ёх°хэш  уЁрэшўэющ чрфрўш ╩ю°ш}
}

\medskip
\noindent
{\it ┬.\,┬. ┴рёют}

\medskip
\noindent
{\footnotesize
╤рэъЄ-╧хЄхЁсєЁуёъшщ уюёєфрЁёЄтхээ√щ єэштхЁёшЄхЄ, \\
╨юёёшщёър  ╘хфхЁрЎш , 199034, ╤рэъЄ-╧хЄхЁсєЁу, ╙эштхЁёшЄхЄёър  эрс., 7-9}

\medskip
{\small \hangindent=\parindent $\quad\ $
╨рёёьрЄЁштрхЄё  юс√ъэютхээюх фшЇЇхЁхэЎшры№эюх єЁртэхэшх яхЁтюую яюЁ фър, ЁрчЁх°хээюх юЄэюёшЄхы№эю яЁюшчтюфэющ. 
╧ЁхфяюырурхЄё , ўЄю хую яЁртр  ўрёЄ№ юяЁхфхыхэр ш эхяЁхЁ√тэр эр ьэюцхёЄтх, 
ёюёЄю ∙хь шч юсырёЄш фтєьхЁэюую хтъышфютр яЁюёЄЁрэёЄтр ш эхъюЄюЁющ ўрёЄш хх уЁрэшЎ√. 
╬с∙хшчтхёЄэю, ўЄю фы  ы■сющ Єюўъш юсырёЄш ш фы  ы■сюую юЄЁхчър ╧хрэю, яюёЄЁюхээюую яю ¤Єющ Єюўъх, 
ЄхюЁхьр ╧хрэю урЁрэЄшЁєхЄ ёє∙хёЄтютрэшх Ёх°хэш  чрфрўш ╩ю°ш эр юЄЁхчъх ╧хрэю.
┬ ёЄрЄ№х ёЇюЁьєышЁютрэ√ єёыютш , яЁш т√яюыэхэшш ъюЄюЁ√ї ьхЄюфюь ыюьрэ√ї ▌щыхЁр фюърчрэю ёє∙хёЄтютрэшх Ёх°хэш  чрфрўш ╩ю°ш, 
яюёЄртыхээющ т уЁрэшўэющ Єюўъх юсырёЄш, ъюЄюЁюх юяЁхфхыхэю эр эхъюЄюЁюь рэрыюух юЄЁхчър ╧хрэю, 
р Єръцх яЁштхфхэ√ єёыютш , урЁрэЄшЁє■∙шх юЄёєЄёЄтшх Єръюую Ёх°хэш . 

{\it ╩ы■ўхт√х ёыютр:} \ уЁрэшўэр  чрфрўр ╩ю°ш, ёє∙хёЄтютрэшх Ёх°хэш , юЄЁхчюъ ╧хрэю. 
}

\bigskip 
{\bf 1. ┬тхфхэшх.} \ 

╨рёёьюЄЁшь юс√ъэютхээюх фшЇЇхЁхэЎшры№эюх єЁртэхэшх яхЁтюую яюЁ фър, ЁрчЁх°хээюх юЄэюёшЄхы№эю яЁюшчтюфэющ
$$y'=f(x,y),\eqno (1)$$
уфх $f(x,y)$ --- тх∙хёЄтхээр  ЇєэъЎш , юяЁхфхыхээр  ш эхяЁхЁ√тэр  эр ьэюцхёЄтх \,$\widetilde G=G\cup \widehat G,$\,
т ъюЄюЁюь $G$ --- ¤Єю юсырёЄ№ т Єюяюыюушш хтъышфютр яЁюёЄЁрэёЄтр $\mathbb{R}^2,$
р ьэюцхёЄтю $\widehat G,$ тючьюцэю яєёЄюх, яЁшэрфыхцшЄ уЁрэшЎх $\partial G$ юсырёЄш $G.$ 
╧Ёш ¤Єюь ъ ьэюцхёЄтє $\widehat G$ сєфхь юЄэюёшЄ№ тёх Єюўъш $\partial G,$ 
т ъюЄюЁ√ї ЇєэъЎш  $f(x,y)$ ьюцхЄ с√Є№ фююяЁхфхыхэр ё ёюїЁрэхэшхь эхяЁхЁ√тэюёЄш. 

┬√сюЁ ьэюцхёЄтр $\widetilde G$ т ърўхёЄтх юсырёЄш юяЁхфхыхэш  ЇєэъЎшш $f(x,y)$ юсєёыютыхэ Єхь, ўЄю яючтюы хЄ 
ЁрёёьрЄЁштрЄ№ т ърўхёЄтх Ёх°хэш  єЁртэхэш  (1) ЇєэъЎш■, уЁрЇшъ ъюЄюЁющ ыхцшЄ тю тёхь ьэюцхёЄтх $\widetilde G,$ 
ш ёЄртшЄ№ чрфрўє ╩ю°ш эх Єюы№ъю т юсырёЄш $G,$ эю ш эр хх уЁрэшЎх. 

╓хы№ яЁхфырурхьющ ЁрсюЄ√ чръы■ўрхЄё  т Єюь, ўЄюс√ т√яшёрЄ№ єёыютш , яЁш ъюЄюЁ√ї ёє∙хёЄтєхЄ Ёх°хэшх чрфрўш ╩ю°ш єЁртэхэш  (1), 
яюёЄртыхээющ т Єюўъх $(x_0,y_0)\in \widehat G,$ ш єёыютш , яЁш ъюЄюЁ√ї Єръюх Ёх°хэш  юЄёєЄёЄтєхЄ, ўЄю, яюэ Єэю, эхтючьюцэю, фы  тэєЄЁхээхщ Єюўъш.

┬ ўрёЄш фюърчрЄхы№ёЄтр ёє∙хёЄтютрэш  яюёЄртыхээр  чрфрўр сєфхЄ Ёх°рЄ№ё  яєЄхь т√фхыхэш  тёхї ёыєўрхт, т ъюЄюЁ√ї тючьюцэю яюёЄЁюшЄ№ рэрыюу юЄЁхчър ╧хрэю ш ёююЄтхЄёЄтє■∙хую ЄЁхєуюы№эшър, ўЄю яючтюышЄ яЁшьхэшЄ№ яЁръЄшўхёъш схч шчьхэхэшщ ьхЄюф ыюьрэ√ї ▌щыхЁр.

╧Ёштхфхь ёэрўрыр эхюсїюфшь√х юяЁхфхыхэш .

\smallskip
{\bf ╬яЁхфхыхэшх 1.}  
╘єэъЎш■ $y=\varphi(x),$ юяЁхфхыхээє■ эр яЁюьхцєЄъх $\langle a,b\,\rangle ,$
сєфхь эрч√трЄ№ {\it Ёх°хэшхь фшЇЇхЁхэЎшры№эюую єЁртэхэш } $(1),$ хёыш т√яюыэ ■Єё  ёыхфє■∙шх єёыютш : \ 
$1)$ ЇєэъЎш  $\varphi(x)$ фшЇЇхЁхэЎшЁєхьр т ы■сющ Єюўъх $x\in\,\langle a,b\,\rangle ;$ 
$2)$ Єюўър $(x,\varphi(x))\in \widetilde G$ яЁш тёхї $x\in\,\langle a,b\,\rangle ;$
$3)$ $\varphi'(x)=f(x,\varphi(x))$ фы  тё ъюую $x\in\,\langle a,b\,\rangle .$

\smallskip
╥хяхЁ№ яюэ Єшх Ёх°хэш  єЁртэхэш  (1) ьюцэю єЄюўэшЄ№ т чртшёшьюёЄш юЄ Єюую, ъръ Ёрёяюыюцхэ хую уЁрЇшъ т ьэюцхёЄтх $\widetilde G.$ 

\smallskip
{\bf ╬яЁхфхыхэшх 2.} \ ╨х°хэшх $y=\varphi(x)$ єЁртэхэш  $(1),$ юяЁхфхыхээюх эр яЁюьхцєЄъх 
$\langle a,b\,\rangle,$ сєфхь эрч√трЄ№:
$a)$ {\it тэєЄЁхээшь Ёх°хэшхь,}\, хёыш Єюўър $(x,\varphi(x))\in G$ фы  ы■сюую $x\in \langle a,b\,\rangle;$ 
$b)$~{\it уЁрэшўэ√ь Ёх°хэшхь,}\, хёыш $(x,\varphi(x))\in \widehat G$ фы  ы■сюую $x\in \langle a,b\,\rangle;$ 
$c)$~{\it ёьх°рээ√ь Ёх°хэшхь,}\, хёыш эрщфєЄё  Єръшх $x_1,x_2\in\langle a,b\,\rangle,$ ўЄю $(x_1,\varphi(x_1))\in G,$ 
р $(x_2,\varphi(x_2))\in \widehat G.$ 
	
\smallskip
╨рчєьххЄё , хёыш єЁртэхэшх (1) ЁрёёьрЄЁштрЄ№ Єюы№ъю т юсырёЄш $G,$ Єю яюэ Єш  Ёх°хэш  ш тэєЄЁхээхую Ёх°хэш  сєфєЄ ёютярфрЄ№.

\smallskip 
{\bf ╟рьхўрэшх 1.} \ 
╘єэъЎшш $f(x,y),$ тёЄЁхўр■∙шхё  яЁш Ёх°хэшш ъюэъЁхЄэ√ї єЁртэхэшщ (1), ъръ яЁртшыю, 
яЁхфёЄрты ■Є ёюсющ ъюьяючшЎш■ ¤ыхьхэЄрЁэ√ї ЇєэъЎшш, эхяЁхЁ√тэ√ї эр ётюхщ юсырёЄш юяЁхфхыхэш , 
ъюЄюЁр  ўр∙х тёхую  ты хЄё  эх сюыхх ўхь ёўхЄэ√ь юс·хфшэхэшхь ёт чэ√ї ьэюцхёЄт $\widetilde G,$ тючьюцэю шьх■∙шї юс∙шх уЁрэшЎ√. 
╚~¤Єю  ты хЄё  юфэющ шч яЁшўшэ ЁрёёьюЄЁхэш  єЁртэхэшщ (1) эх т юсырёЄш, р эр ьэюцхёЄтх $\widetilde G.$
═ряЁшьхЁ, т єЁртэхэшш $y'=\sqrt{y},$ юсырёЄ№■ юяЁхфхыхэш  ъюЄюЁюую  ты хЄё  ьэюцхёЄтю $\widetilde G=\{(x,y)\!:\,x\in \mathbb{R}^1,\,y\ge 0\},$  ёюуырёэю яЁштхфхээ√ь т√°х юяЁхфхыхэш ь ЇєэъЎш  $y(x)\equiv 0$ эр $\mathbb{R}^1$  ты хЄё  Ёх°хэшхь шыш уЁрэшўэ√ь Ёх°хэшхь. 

\smallskip
{\bf ╬яЁхфхыхэшх 3.} \ 
╟рфрўє ╩ю°ш ё эрўры№э√ьш фрээ√ьш (э.\,ф.) $x_0,y_0$ сєфхь эрч√трЄ№ {\it тэєЄЁхээхщ чрфрўхщ ╩ю°ш,}\, хёыш Єюўър $(x_0,y_0)\in G,$ 
ш сєфхь эрч√трЄ№ {\it уЁрэшўэющ чрфрўхщ ╩ю°ш,}\, хёыш Єюўър $(x_0,y_0)\in \widehat G.$ 

\smallskip
{\bf ╬яЁхфхыхэшх 4.}  {\it ┬эєЄЁхээхх (уЁрэшўэюх, ёьх°рээюх) Ёх°хэшх чрфрўш ╩ю°ш} єЁртэхэш  $(1.1)$ ё эрўры№э√ьш фрээ√ьш 
$x_0,y_0$ {\it ёє∙хёЄтєхЄ,}\, хёыш Єюўър $(x_0,y_0)\in G\,(\widehat G,\widetilde G)$ ш эрщфєЄё  яЁюьхцєЄюъ $\langle a,b\,\rangle,$ 
ёюфхЁцр∙шщ $x_0,$ ш юяЁхфхыхээюх эр эхь тэєЄЁхээхх (уЁрэшўэюх, ёьх°рээюх) Ёх°хэшх $y=\varphi(x)$ Єръшх, ўЄю $\varphi(x_0)=y_0.$  

\smallskip
╥ръшь юсЁрчюь, уЁрЇшъ тэєЄЁхээхую Ёх°хэш  чрфрўш ╩ю°ш ыхцшЄ т юсырёЄш $G,$ уЁрэшўэюую --- т $\widehat G,$ 
р ёьх°рээюую --- ш т юсырёЄш $G,$ ш т ьэюцхёЄтх $\widehat G.$

╧хЁхїюф  ъ тюяЁюёє ю ёє∙хёЄтютрэшш Ёх°хэш  чрфрўш ╩ю°ш, яЁштхфхь ЇюЁьєышЁютъє їюЁю°ю шчтхёЄэющ ЄхюЁхь√ ╧хрэю, 
тёЄЁхўр■∙є■ё  яЁръЄшўхёъш тю тёхї ъэшурї яю юс√ъэютхээ√ь фшЇЇхЁхэЎшры№э√ь єЁртэхэш ь (ёь.\,эряЁ.\,[1]--[5]) 
ш ърёр■∙є■ё  т яЁхфыюцхээющ ЄхЁьшэюыюушш ёє∙хёЄтютрэш  Ёх°хэш  тэєЄЁхээхщ чрфрўш ╩ю°ш. 

\smallskip 
{\bf ╥хюЁхьр ╧хрэю.} \ 
{\it ╧єёЄ№ т єЁртэхэшш (1) ЇєэъЎш  $f(x,y)$ 
юяЁхфхыхэр ш эхяЁхЁ√тэр т юсырёЄш $G,$ Єюуфр фы  ы■сющ Єюўъш $(x_0,y_0)\in G$ 
ш фы  ы■сюую юЄЁхчър ╧хрэю $P_h(x_0,y_0)$ ёє∙хёЄтєхЄ яю ъЁрщэхщ ьхЁх юфэю Ёх°хэшх чрфрўш ╩ю°ш єЁртэхэш  $(1)$ 
ё эрўры№э√ьш фрээ√ьш $x_0,y_0,$ юяЁхфхыхээюх эр $P_h(x_0,y_0).$ }

{\small ╟фхё№ юЄЁхчюъ ╧хрэю $P_h(x_0,y_0)=[x_0-h,x_0+h],$ р ъюэёЄрэЄр $h>0$ юяЁхфхы хЄё  ёыхфє■∙шь юсЁрчюь: ёє∙хёЄтє■Є ъюэёЄрэЄ√ 
$a,b>0$ Єръшх, ўЄю яЁ ьюєуюы№эшъ $R=\{(x,y)\,:\,|x-x_0|\le a,\,|y-y_0|\le b\}\subset G.$
┼ёыш $f(x,y)\equiv 0$ эр $\overline R,$ Єю $h=a.$
┬ яЁюЄштэюь ёыєўрх яюыюцшь \,$M=\max_{(x,y)\in R}|f(x,y)|>0,$ Єюуфр $h=\min\,\{a,b/M\}.$ 
├хюьхЄЁшўхёъш $h$ --- фышэр т√ёюЄ√ ЄЁхєуюы№эшър $T^+$ ё тхЁ°шэющ т Єюўъх $(x_0,y_0),$ сюъют√ьш ёЄюЁюэрьш, 
шьх■∙шьш Єрэухэё√ єуыют эръыюэр $\pm M,$ ш юёэютрэшхь, ыхцр∙хь эр яЁ ьющ $x=x_0+h$ $(T^-$ ёЄЁюшЄё  рэрыюушўэю). }

\smallskip 
{\bf ╟рьхўрэшх 2.} \ ┬ [6] яЁхфяюырурхЄё , ўЄю ЇєэъЎш  $f(x,y)$ юяЁхфхыхэр ш эхяЁхЁ√тэр т яЁ ьюєуюы№эшъх 
$R^+=\{(x,y)\,:\,x_0\le x\le x_0+a,\,|y-y_0|\le b\},$ ўЄю, эр ёрьюь фхых, 
Єюы№ъю ш ЄЁхсєхЄё  яЁш ы■сюь трЁшрэЄх фюърчрЄхы№ёЄтр ЄхюЁхь√ ╧хрэю фы  $x\ge x_0.$

\smallskip
┬ ЁрсюЄх сєфхЄ ЁрёёьрЄЁштрЄ№ё  чрфрўр ╩ю°ш, яюёЄртыхээр  т яЁюшчтюы№эющ уЁрэшўэющ Єюўъх $(x_0,y_0)$ ьэюцхёЄтр $\widetilde G.$
═ю яЁхцфх, ўхь ёфхырЄ№ яЁхфяюыюцхэш , ърёр■∙шхё  ёЄЁєъЄєЁ√ уЁрэшЎ√ т ьрыющ юъЁхёЄэюёЄш ¤Єющ Єюўъш, 
фы  єяЁю∙хэш  шёяюы№чєхь√ї т фры№эхщ°хь юсючэрўхэшщ ш ЇюЁьєы сєфхь ёўшЄрЄ№, эх єьхэ№°р  юс∙эюёЄш,
ўЄю чрфрўр ╩ю°ш яюёЄртыхэр 
т Єюўъх $(0,0)$ ш ўЄю ЇєэъЎш  $f$ т эхщ Ёртэр эєы■, Є.\,х. сєфхь ЁрёёьрЄЁштрЄ№ єЁртэхэшх
$$y'=f_0(x,y),\eqno (2)$$ 
т ъюЄюЁюь ЇєэъЎш  $f_0$ юяЁхфхыхэр ш эхяЁхЁ√тэр эр ьэюцхёЄтх $\widetilde G=G\cup \widehat G,$ 
$G$ --- юсырёЄ№ т $\mathbb{R}^2,$ $\widehat G\in \partial G,$ 
Єюўър $\text{O}=(0,0)\in \widehat G,$ $f_0(0,0)=0$ ш яюёЄртыхэр чрфрўр ╩ю°ш ё эрўры№э√ьш фрээ√ьш $0,0.$ 

\smallskip
{\small ┬ ёрьюь фхых, тюч№ьхь яЁюшчтюы№эюх єЁртэхэшх (1) ш яюёЄртшь чрфрўє ╩ю°ш т ы■сющ Єюўъх $(x_0,y_0)\in \widehat G,$ 
Єюуфр чрьхэр \,$x=t+x_0,\quad y=v+y_0+f(x_0,y_0)t$\, ётхфхЄ єЁртэхэшх (1) ъ єЁртэхэш■ \,$\dot v=f_0(t,v),$\, 
т ъюЄюЁюь ЇєэъЎш  $f_0(t,v)=f(t+x_0,v+y_0+f(x_0,y_0)t)-f(x_0,y_0).$
╤ыхфютрЄхы№эю, яЁш $x=x_0,\ y=y_0$ яюыєўрхь $t=t_0=0,\ v=v_0=0$ ш $f_0(0,0)=0.$ }

\medskip
{\bf 2. ├Ёрэшўэ√х ъЁшт√х ш яюЁюцфрхь√х шьш ьэюцхёЄтр.} 

\smallskip
{\bf ╬яЁхфхыхэшх 5.} \ ╘єэъЎш■, юяЁхфхыхээє■ эр эхъюЄюЁюь юЄЁхчъх $[0,a_u],$ сєфхь эрч√трЄ№ {\it яЁртющ тхЁїэхуЁрэшўэющ ЇєэъЎшхщ}\, 
ш юсючэрўрЄ№ $y=b_{a_u}^+(x),$ 
хёыш т√яюыэ ■Єё  ёыхфє■∙шх я Є№ єёыютшщ:   
\,$1)\ b_{a_u}^+(x)\in C^1([0,a_u]);$ 
$2)\ b_{a_u}^+(0)=0;$ 
$3)$~${b_{a_u}^+}'(0)\ge 0;$ 
$4)$~$b_{a_u}^+$ т√яєъыр тэшч эр $[0,a_u],$ хёыш ${b_{a_u}^+}'(0)=0;$ 
$5)$~уЁрЇшъ $b_{a_u}^+$ --- {\it яЁртр  тхЁїэхуЁрэшўэр  ъЁштр }    
$\gamma_{a_u}^+=\{x\in [0,a_u],\,y=b_{a_u}^+(x)\}$ --- яЁшэрфыхцшЄ ьэюцхёЄтє $\widehat G.$ 

{\small ╟фхё№ ёшьтюы \,$+$\, яюфЁрчєьхтрхЄ {\sl яЁртюёЄюЁюээшщ}, $b$ --- {\sl уЁрэшўэр  ЇєэъЎш }, $u$ --- {\sl тхЁїэшщ}, $l$ --- {\sl эшцэшщ}, 
р т√яєъыюёЄ№ яюэшьрхЄё  т эхёЄЁюуюь ёь√ёых, Є.\,х. фюяєёър■Єё  ЄюцфхёЄтр $b_{a_u}^+(x)\equiv 0$ эр $[0,a_u]$ шыш $b_{a_l}^+(x)\equiv 0$ эр $[0,a_l].$}

┬ Ёхчєы№ЄрЄх $\gamma_{a_u}^+$ --- уырфър  ъЁштр  шч $\widehat G,$ ярЁрьхЄЁшчютрээр  ЇєэъЎшхщ $b_{a_u}^+(x).$ 
╬эр эрўшэрхЄё  т Єюўъх $\text{O}$ ш Ёрёяюыюцхэр т яхЁтющ ўхЄтхЁЄш, р ъЁштр  $\gamma_{a_l}^+,$ эрўшэр ё№ Єрь цх, Ёрёяюыюцхэр т ўхЄтхЁЄющ ўхЄтхЁЄш.

└эрыюушўэю ттюф Єё  {\it яЁрт√х эшцэхуЁрэшўэр  ЇєэъЎш }\, $y=b_{a_l}^+(x)$ ш {\it яЁртр  эшцэхуЁрэшўэр  ъЁштр } $\gamma_{a_l}^+,$ 
Єюы№ъю т $3)$ ${b_{a_l}^+}'(0)\le 0,$ р т $4)$ $b_{a_l}^+(x)$ т√яєъыр ттхЁї. 

\smallskip
═хяюёЁхфёЄтхээю шч юяЁхфхыхэш  т√ЄхърхЄ, ўЄю фы  ы■сющ яЁртющ тхЁїэхуЁрэшўэющ ЇєэъЎшш $b_{a_u^*}^+(x)$ 
ЇєэъЎш  $b_{a_u}^+(x),$  ты ■∙р ё  хх ёєцхэшхь эр яЁюшчтюы№э√щ юЄЁхчюъ $[0,a_u]$ ё $a_u<a_u^*,$ 
юёЄрхЄё  яЁртющ тхЁїэхуЁрэшўэющ. ─ы  $b_{a_l^*}^+(x)$ --- тёх рэрыюушўэю. 

╧юыюцшь $\tau_u={b_{a_u}^+}'(0)/2,\ \tau_l=-{b_{a_l}^+}'(0)/2,$ ш, эх єьхэ№°р  юс∙эюёЄш, сєфхь ёўшЄрЄ№, ўЄю т√яюыэ ■Єё  єёыютш : 
$$\begin{matrix} 
b_{a_u}^+(a_u)\le a_u\ \text{ яЁш }\ \tau_u=0,\ \   
\forall\,x\in [0,a_u]\!:\ {b_{a_u}^+}'(x)\ge \tau_u\ \text{ яЁш }\ \tau_u>0; \\
-b_{a_l}^+(a_u)\le a_l\ \text{ яЁш }\ \tau_l=0,\ \   
\forall\,x\in [0,a_l]\!:\ \;{b_{a_l}^+}'(x)\le \tau_l\ \text{ яЁш }\ \tau_l>0.\hfill 
  \end{matrix}\eqno (3^+)$$ 

─ы  тё ъюую $c>0$ ЁрёёьюЄЁшь яЁртє■ $c$-юъЁхёЄэюёЄ№ Єюўъш $\text{O}:$ 
$$N_c^+=\{(x,y)\!:\,x\in(0,c],\,|y|\le c\}.$$ 

╬ЄьхЄшь ёЁрчє, ўЄю эхЁртхэёЄтю $b_{a_u}^+(a_u)\le a_u$ шыш $-b_{a_l}^+(a_l)\le a_l$ шч $(3^+)$ яЁш тёхї $c\le a_u$ шыш $c\le a_l$ урЁрэЄшЁєхЄ яхЁхёхўхэшх ъЁштющ $\gamma_{a_u}^+$ шыш $\gamma_{a_l}^+$ ё сюъютющ, р эх эр тхЁїэхщ шыш эшцэхщ ёЄюЁюэющ яЁ ьюєуюы№эшър $N_c^+.$ 

\smallskip 
┬тхфхь ЄЁш Єшяр ьэюцхёЄт, ъюЄюЁ√х юяЁхфхы ■Єё  яЁрт√ьш тхЁїэхуЁрэшўэ√ьш ш эшцэхуЁрэшўэ√ьш ъЁшт√ьш.

1. ─ы  ы■сющ яЁртющ тхЁїэхуЁрэшўэющ ъЁштющ $\gamma_{a_u^*}^+,$ ярЁрьхЄЁшчєхьющ ЇєэъЎшхщ $b_{a_u^*}^+(x),$ 
фы  ъюЄюЁющ т√яюыэ хЄё  єёыютшх$(3_u^+),$ яюыюцшь $c_u^*=\max\{a_u^*,b_{a_u^*}^+(a_u^*)\}.$ ╥юуфр 
$$\forall\,c\ (0<c\le c_u^*)\quad \exists\,a_u\ (0<a_u\le a_u^*)\!:\ \ \max\{a_u,b_{a_u}^+(a_u)\}=c,$$ 
уфх ЇєэъЎш  $b_{a_u}^+(x)$ --- ¤Єю ёєцхэшх $b_{a_u^*}^+(x).$
╥ръшь юсЁрчюь, Єюўър $(a_u,b_{a_u}^+(a_u))$ --- ¤Єю яЁрт√щ ъюэхЎ тхЁїэхуЁрэшўэющ ъЁштющ $\gamma_{a_u}^+.$ 
╬э Ёрёяюыюцхэ эр тхЁїэхщ ёЄюЁюэх яЁ ьюєуюы№эшър $N_c^+,$ хёыш $c=b_{a_u}^+(a_u)\ge a_u,$ ш эр сюъютющ, хёыш $b_{a_u}^+(a_u)\le a_u=c.$

╥хяхЁ№ фы  тё ъюую $c\ \ (0<c\le c_u^*)$ ттхфхь ьэюцхёЄтю 
$$U_c^+=\{(x,y)\!:\, x\in (0,a_u],\ -c\le y\le b_{a_u}^+(x);\ x\in (a_u,c],\ |y|\le c\}$$
--- ¤Єю яЁ ьюєуюы№эшъ $N_c^+,$ шч ъюЄюЁюую "т√Ёхчрэ"\ эрфуЁрЇшъ ъЁштющ $\gamma_{a_u}^+.$ 

2. ─ы  ы■сющ яЁртющ эшцэхуЁрэшўэющ ъЁштющ $\gamma_{a_l^*}^+,$ ярЁрьхЄЁшчєхьющ ЇєэъЎшхщ $b_{a_l^*}^+(x),$ 
фы  ъюЄюЁющ т√яюыэ хЄё  єёыютшх$(3_l^+),$ яюыюцшь$c_l^*=\max\{a_l^*,-b_{a_l^*}^+(a_l^*)\}.$ ╥юуфр
$$\forall\,c\ (0<c\le c_l^*)\quad \exists\,a_l\ (0<a_l\le a_l^*)\!:\ \ \max\{a_l,-b_{a_l}^+(a_l)\}=c,$$ 
уфх ЇєэъЎш  $b_{a_l}^+(x)$ --- ¤Єю ёєцхэшх $b_{a_l^*}^+(x).$
╥ръшь юсЁрчюь, Єюўър $(a_l,b_{a_l}^+(a_l))$ --- ¤Єю яЁрт√щ ъюэхЎ ъЁштющ $\gamma_{a_l}^+.$ 
╬э Ёрёяюыюцхэ эр эшцэхщ ёЄюЁюэх яЁ ьюєуюы№эшър $N_c^+,$ хёыш $c=-b_{a_l}^+(a_l)\ge a_l,$ ш эр сюъютющ, хёыш $-b_{a_l}^+(a_l)\le a_l=c.$

╥хяхЁ№ фы  тё ъюую $c\ \ (0<c\le c_l^*)$ ттхфхь ьэюцхёЄтю 
$$O_c^+=\{(x,y)\!:\, x\in (0,a_l],\ b_{a_l}^+(x)\le y\le c;\ x\in (a_l,c],\ |y|\le c\}$$
--- ¤Єю яЁ ьюєуюы№эшъ $N_c^+,$ шч ъюЄюЁюую "т√Ёхчрэ"\ яюфуЁрЇшъ ъЁштющ $\gamma_{a_l}^+.$ 

3. ─ы  ы■с√ї яЁртющ тхЁїэхуЁрэшўэющ ъЁштющ $\gamma_{a_u^*}^+,$ яЁртющ эшцэхуЁрэшўэющ ъЁштющ $\gamma_{a_l^*}^+$ 
ш тё ъюую $c\ (0<c\le c_b^*),$ уфх $c_b^*=\min\{c_u^*,c_l^*\},$ р ъюэёЄрэЄ√ $c_u^*,c_l^*$ юяЁхфхыхэ√ т√°х, ттхфхь ьэюцхёЄтю 
$$B_c^+=U_c^+\cap O_c^+$$ 
--- ¤Єю $N_c^+,$ шч ъюЄюЁюую "т√Ёхчрэ√"\ эрфуЁрЇшъ ъЁштющ $\gamma_{a_u}^+$ ш яюфуЁрЇшъ ъЁштющ $\gamma_{a_l}^+.$ 

{\small ╟фхё№ $U$ яюфЁрчєьхтрхЄ {\sl яюф}, $O$ --- {\sl эрф}, $B$ --- {\sl ьхцфє}. 

\smallskip
{\bf ╬яЁхфхыхэшх 6.} \ ┴єфхь уютюЁшЄ№, ўЄю фы  єЁртэхэш  (2) ЁхрышчєхЄё : 
 
{\it ёыєўрщ} $U^+],$ хёыш $\exists\, c_u\ (0<c_u \le c_u^*)\!:\ U_{c_u}^+\cap \widehat G=\gamma_{a_u}^+\backslash \text{O};$

{\it ёыєўрщ} $O^+],$ хёыш $\exists\, c_l\ (0<c_l \le c_l^*)\!:\ O_{c_l}^+\cap \widehat G=\gamma_{a_l}^+\backslash \text{O};$

{\it ёыєўрщ} $B^+],$ хёыш $\exists\, c_b\ (0<c_b \le c_b^*)\!:\ 
  B_{c_b}^+\cap \widehat G=(\gamma_{a_u}^+\cup \gamma_{a_l}^+)\backslash \text{O}.$
	
\smallskip
┬ фры№эхщ°хь т юсючэрўхэш ї ьэюцхёЄтр $X_c^+$ ш ёыєўр  $X^+]$ сєътр $X$ сєфхЄ яюфЁрчєьхтрЄ№ ы■сє■ шч сєът $U,\,O$ шыш $B.$

┬ ёыєўрх $X^+]$ фы  ьэюцхёЄтр $X_{c_*}^+$ ышсю $X_{c_*}^+\cap G\ne \emptyset,$ 
ўЄю Ёртэюёшы№эю Єюьє, ўЄю $X_{c_*}^+$ схч тїюф ∙шї т эхую уЁрэшўэ√ї ъЁшт√ї ыхцшЄ т $G,$ ышсю $X_{c_*}^+\cap G=\emptyset.$  

{\small ╟фхё№ ёшьтюы $*$ --- ¤Єю ёююЄтхЄёЄтє■∙р  $X$ сєътр $u,\,l$ шыш $b.$}

┬ ёт чш ё ¤Єшь ёыєўрщ $X^+]$ т чртшёшьюёЄш юЄ Ёрёяюыюцхэш  ьэюцхёЄтр $X_{c_*}^+$ юЄэюёшЄхы№эю юсырёЄш $G$ 
ЁрёярфрхЄё  эр фтр ёыєўр , ъюЄюЁ√х сєфхь юсючэрўрЄ№ $X_1^+]$ шыш $X_2^+].$ 
╧Ёш ¤Єюь фюяюыэшЄхы№э√щ шэфхъё $>,\,=$ шыш $<$ яЁш хую эрышўшш т юсючэрўхэшш ы■сюую шч °хёЄш тючэшъ°шї ёыєўрхт 
сєфхЄ єЄюўэ Є№ чэръ яЁюшчтюфэющ ёююЄтхЄёЄтє■∙хщ тхЁїэхуЁрэшўэющ шыш (ш) эшцэхуЁрэшўэющ ЇєэъЎшш т эєых. 

─рфшь ёЄЁюуюх юяшёрэшх тючэшър■∙шї ёыєўрхт, ёюяЁютюфшт тёх ЁрчэютшфэюёЄш ёыєўрхт $X_1^+]$ яю ёэ ■∙шьш Ёшёєэърьш.

$U_1^+]:$ $(U_{c_u}^+\backslash \gamma_{a_u}^+)\subset G,$\,  
2 яюфёыєўр : $U_1^{+,>}]:\ {b_{a_u}^+}'(0)>0,\ \ U_1^{+,=}]:\ {b_{a_u}^+}'(0)=0;$ 

$U_2^+]:$ $U_{c_u}^+\cap G=\emptyset,$ яюфёыєўрш Єх цх. 

$O_1^+]:$ $(O_{c_l}^+\backslash \gamma_{a_l}^+)\subset G,$\,  
2 яюфёыєўр : $O_{1,<}^+]:\ {b_{a_l}^+}'(0)<0,\ \ O_{1,=}^+]:\ {b_{a_l}^+}'(0)=0;$ 

$O_2^+]:$ $O_{c_l}^+\cap G=\emptyset,$ яюфёыєўрш Єх цх. 

$B_1^+]:$ $(B_{c_b}^+\backslash (\gamma_{a_u}^+ \cup \gamma_{a_l}^+))\subset G,$\,  
4 яюфёыєўр : 
$B_{1,<}^{+,>}]:\ {b_{a_u}^+}'(0)>0,\,{b_{a_l}^+}'(0)<0,$  
$B_{1,=}^{+,=}]:\ {b_{a_u}^+}'(0)=0,\,{b_{a_l}^+}'(0)=0,$ 
$B_{1,=}^{+,>}]:\ {b_{a_u}^+}'(0)>0,\,{b_{a_l}^+}'(0)=0,$ 
$B_{1,<}^{+,=}]:\ {b_{a_u}^+}'(0)=0,\,{b_{a_l}^+}'(0)<0;$ 

$B_2^+]:$ $B_{c_b}^+\cap G=\emptyset,$ яюфёыєўрш Єх цх. 

\begin{center}
\begin{tikzpicture}[scale=1.]
\tikzset{line01/.style={line width=0.7pt}}
\tikzset{line03/.style={line width=1pt}}
\useasboundingbox (0.7,-2) rectangle (2.9,2.6);
\fill[whitesmoke] (0,0) ..controls +(60:0.5cm) and +(210:0.6cm).. (1.5,2) -- (2,2) -- (2,-2) -- (0,-2) -- cycle;
\draw[->] (-0.3,0) -- (2.8,0) node[below] {\small $x$};
\draw[->] (0,-2.2) -- (0,2.4) node[left] {\small $y$};
\draw (0,-2) node [above right] {$U_{c_u}^+$} rectangle (2,2);
\draw[dashed] (1.5,0) node[below=-2pt] {\small $a_u$} -- (1.5,2); 
\draw[line01] (0,0) ..controls +(60:0.5cm) and +(210:0.6cm).. (1.5,2) node[near end,left=1pt] {$\gamma_{a_u}^+$};
\draw (0,0)+(-135:7pt) node {\small 0};
\draw (2,0)+(-45:7pt) node {\small $c_u$};
\draw (0,2) node [left] {\small $c_u$};
\draw (0,-2) node [left] {\small $-c_u$};
\draw (1.0,2.3) node {$U_1^{+,>}]$};
\end{tikzpicture}
\begin{tikzpicture}[scale=1.]
\tikzset{line01/.style={line width=0.7pt}}
\tikzset{line03/.style={line width=1pt}}
\useasboundingbox (-0.4,-2) rectangle (2.9,2.6);
\fill[whitesmoke] (0,0) ..controls +(0:0.5cm) and +(250:0.6cm).. (2,1.6) -- (2,-2) -- (0,-2) -- cycle;
\draw[->] (-0.3,0) -- (2.8,0) node[below] {\small $x$};
\draw[->] (0,-2.2) -- (0,2.4) node[left] {\small $y$};
\draw (0,-2) node [above right] {$U_{c_u}^+$} rectangle (2,2);
\draw[line01] (0,0) ..controls +(0:0.5cm) and +(250:0.6cm).. (2,1.6) node[near end,left=2pt] {$\gamma_{a_u}^+$};
\draw (0,0)+(-135:7pt) node {\small 0};
\draw (1.8,-0.17) node {\small $a_u=c_u$};
\draw (0,2) node [left] {\small $c_u$};
\draw (0,-2) node [left] {\small $-c_u$};
\draw (1.0,2.3) node {$U_1^{+,=}]$};
\end{tikzpicture}
\begin{tikzpicture}[scale=1.]
\tikzset{line01/.style={line width=0.7pt}}
\tikzset{line03/.style={line width=1pt}}
\useasboundingbox (-0.4,-2) rectangle (2.9,2.6);
\fill[whitesmoke] (0,0) ..controls +(-60:0.5cm) and +(160:0.6cm).. (2,-1.7) -- (2,2) -- (0,2) -- cycle;
\draw[->] (-0.3,0) -- (2.8,0) node[below] {\small $x$};
\draw[->] (0,-2.2) -- (0,2.4) node[left] {\small $y$};
\draw (0,-2) rectangle (2,2);
\draw (0,2) node [below right] {$O_{c_l}^+$};
\draw[line01] (0,0) ..controls +(-60:0.5cm) and +(160:0.6cm).. (2,-1.7) node[near end,left=6pt] {$\gamma_{a_l}^+$};
\draw (0,0)+(-135:7pt) node {\small 0};
\draw (1.8,-0.17) node {\small $a_l=c_l$};
\draw (0,2) node [left] {\small $c_l$};
\draw (0,-2) node [left] {\small $-c_l$};
\draw (1.0,2.3) node {$O_{1,<}^+]$};
\end{tikzpicture}
\begin{tikzpicture}[scale=1.]
\tikzset{line01/.style={line width=0.7pt}}
\tikzset{line03/.style={line width=1pt}}
\useasboundingbox (-0.4,-2) rectangle (1.7,2.6);
\fill[whitesmoke] (0,0) ..controls +(0:0.5cm) and +(110:0.6cm).. (2,-1.6) -- (2,2) -- (0,2) -- cycle;
\draw[->] (-0.3,0) -- (2.8,0) node[below] {\small $x$};
\draw[->] (0,-2.2) -- (0,2.4) node[left] {\small $y$};
\draw (0,-2) rectangle (2,2);
\draw (0,2) node [below right] {$O_{c_l}^+$};
\draw[line01] (0,0) ..controls +(0:0.5cm) and +(110:0.6cm).. (2,-1.6) node[pos=0.9, left=2pt] {$\gamma_{a_l}^+$};
\draw (0,0)+(-135:7pt) node {\small 0};
\draw (1.8,-0.17) node {\small $a_l=c_l$};
\draw (0,2) node [left] {\small $c_l$};
\draw (0,-2) node [left] {\small $-c_l$};
\draw (1.0,2.3) node {$O_{1,=}^+]$};
\end{tikzpicture}
\end{center}		

\begin{center}
\begin{tikzpicture}[scale=1.]
\tikzset{line01/.style={line width=0.7pt}}
\tikzset{line03/.style={line width=1pt}}
\useasboundingbox (0.7,-2) rectangle (2.9,2.6);
\fill[whitesmoke] (0,0) ..controls +(70:0.5cm) and +(210:0.6cm).. (1.5,2) -- (2,2) -- (2,-2) -- (1.8,-2) 
                        ..controls +(100:0.2cm) and +(-70:0.2cm).. 	(1.5,-1.6) 
												..controls +(110:1.4cm) and +(-40:1cm).. cycle;
\draw [->] (-0.3,0) -- (2.8,0) node[below] {\small $x$};
\draw [->] (0,-2.2) -- (0,2.4) node[left] {\small $y$};
\draw (0,-2) rectangle (2,2);
\draw[dashed] (1.5,0) node[below=-2pt] {\small $a_u$} -- (1.5,2); 
\draw[dashed] (1.8,0) node[above=-2pt] {\small $a_l$} -- (1.8,-2); 
\draw[line01] (0,0) ..controls +(70:0.5cm) and +(210:0.6cm).. (1.5,2) node[near end,left=0pt] {$\gamma_{a_u}^+$}
		node [pos=0.25, right=-2pt] {$B_{c_b}^+$};
\draw[line01] (0,0) ..controls +(-40:1cm) and +(110:1.4cm).. (1.5,-1.6) node[pos=0.83, left=-1pt] {$\gamma_{a_l}^+$} 
		                ..controls +(-70:0.2cm) and +(100:0.2cm).. (1.8,-2);
\draw (0,0)+(-135:7pt) node {\small 0};
\draw (2,0)+(-45:7pt) node {\small $c_b$};
\draw (0,2) node [left] {\small $c_b$};
\draw (0,-2) node [left] {\small $-c_b$};
\draw (1.0,2.3) node {$B_{1,<}^{+,>}]$};
\end{tikzpicture}
\begin{tikzpicture}[scale=1.]
\tikzset{line01/.style={line width =0.7pt}}
\tikzset{line03/.style={line width =1pt}}
\useasboundingbox (-0.4,-2) rectangle (2.9,2.6);
\fill[whitesmoke] (0,0) ..controls +(0:0.4cm) and +(250:0.6cm).. (2,1.8) -- (2,-1.7) ..controls +(100:1cm) and +(0:1.5cm).. cycle;
\draw[->] (-0.3,0) -- (2.8,0) node[below] {\small $x$};
\draw[->] (0,-2.2) -- (0,2.4) node[left] {\small $y$};
\draw (0,-2) rectangle (2,2); 
\draw (0,2); 
\draw[line01] (0,0) ..controls +(0:0.4cm) and +(250:0.6cm).. (2,1.8) 
     node[near end, left=+2 pt] {$\gamma_{a_u}^+$} node[pos=0.45, right=0pt] {$B_{c_b}^+$}; 
\draw[line01] (2,-1.7) ..controls +(100:1cm) and +(0:1.5cm).. node[near start, left=-2pt] {$\gamma_{a_l}^+$} (0,0);
\draw (0,0)+(-135:7pt) node {\small 0};
\draw (1.75,0.15) node {\small $a_u=a_l=c_b$};
\draw (0,2) node [left] {\small $c_b$};
\draw (0,-2) node [left] {\small $-c_b$};
\draw (1.0,2.3) node {$B_{1,=}^{+,=}]$};
\end{tikzpicture}
\begin{tikzpicture}[scale=1.]
\tikzset{line01/.style={line width =0.7pt}}
\tikzset{line03/.style={line width =1pt}}
\useasboundingbox (-0.4,-2) rectangle (2.9,2.6);
\fill[whitesmoke] (0,0) ..controls +(60:0.5cm) and +(210:0.6cm).. (1.5,2) -- (2,2) -- (2,-1.75) 
                        ..controls +(100:0.4cm) and +(-60:0.1cm).. (1.5,-0.75) ..controls +(120:0.5cm) and +(0:0.5cm).. (0,0);
\draw[->] (-0.3,0) -- (2.8,0) node[below] {\small $x$};
\draw[->] (0,-2.2) -- (0,2.4) node[left] {\small $y$};
\draw (0,-2) rectangle (2,2); 
\draw (0,2); 
\draw[line01] (1.5,-0.75) ..controls +(120:0.5cm) and +(0:0.5cm).. (0,0);
\draw[dashed] (1.5,0) node[below=-2pt] {\small $a_u$} -- (1.5,2); 
\draw[line01] (0,0) ..controls +(60:0.5cm) and +(210:0.6cm).. (1.5,2) 
    node[near end,left=0pt] {$\gamma_{a_u}^+$} node [pos=0.2, right=0pt] {$B_{c_b}^+$};
\draw[line01] (1.5,-0.75) ..controls +(-60:0.1cm) and +(100:0.4cm).. (2,-1.75) node[pos=0.85, left=-1pt] {$\gamma_{a_l}^+$};
\draw (0,0)+(-135:7pt) node {\small 0};
\draw (2.05,0.15) node {\small $a_l=c_b$};
\draw (0,2) node [left] {\small $c_b$};
\draw (0,-2) node [left] {\small $-c_b$};
\draw (1.0,2.3) node {$B_{1,=}^{+,>}]$};
\end{tikzpicture}
\begin{tikzpicture}[scale=1.]
\tikzset{line01/.style={line width =0.7pt}}
\tikzset{line03/.style={line width =1pt}}
\useasboundingbox (-0.4,-2) rectangle (1.7,2.6);
\fill[whitesmoke] (0,0) ..controls +(-60:0.5cm) and +(160:0.6cm).. (2,-1.7) -- (2,1.6) ..controls +(-120:0.5cm) and +(0:0.5cm).. cycle; 
\draw[->] (-0.3,0) -- (2.8,0) node[below] {\small $x$};
\draw[->] (0,-2.2) -- (0,2.4) node[left] {\small $y$};
\draw (0,-2) rectangle (2,2); 
\draw (0,2); 
\draw[line01] (2,1.6) ..controls +(-120:0.5cm) and +(0:0.5cm)..   node[near start,left=+2pt] {$\gamma_{a_u}^+$} (0,0);
\draw[line01] (0,0) ..controls +(-60:0.5cm) and +(160:0.6cm).. (2,-1.7) 
   node[near end,left=5pt] {$\gamma_{a_l}^+$};
\draw (1.4,-0.5) node {$B_{c_b}^+$}; 
\draw (0,0)+(-135:7pt) node {\small 0};
\draw (1.75,0.15) node {\small $a_u=a_l=c_b$};
\draw (0,2) node [left] {\small $c_b$};
\draw (0,-2) node [left] {\small $-c_b$};
\draw (1.0,2.3) node {$B_{1,<}^{+,=}]$};
\end{tikzpicture}
\end{center}

\smallskip
{\bf ╟рьхўрэшх 3.} \ ┬ ЁрёёьюЄЁхээ√ї ёыєўр ї яЁртр  $c$"=юъЁхёЄэюёЄ№ $N_c^+,$ ъюэхўэю, ьюцхЄ ёюфхЁцрЄ№ сюыхх юфэющ эшцэхуЁрэшўэющ 
ш сюыхх юфэющ тхЁїэхуЁрэшўэющ ъЁштющ. ╬фэръю, эрышўшх шыш юЄёєЄёЄтшх уЁрэшўэ√ї ъЁшт√ї т яЁ ьюєуюы№эшъх $N_c^+$ тэх ьэюцхёЄтр $X_c^+$ 
эх тыш хЄ эр ёє∙хёЄтютрэшх Ёх°хэш  уЁрэшўэющ чрфрўш ╩ю°ш т ёыєўр ї $X_1^+].$ 
┬ Єю цх тЁхь  т ёыєўрх $U_2^{+,=}$ эхёюьэхээ√щ шэЄхЁхё яЁхфёЄрты ■Є ёшЄєрЎшш, ъюуфр т юсырёЄш $G$ ёюфхЁцшЄё  
ъръ тхё№ эрфуЁрЇшъ ъЁштющ $\gamma_{a_u}^+,$ ыхцр∙шщ т $N_{c_u}^+,$ Єръ ш ьэюцхёЄтю, 
чръы■ўхээюх ьхцфє $\gamma_{a_u}^+$ ш х∙х юфэющ яЁртющ тхЁїэхуЁрэшўэющ ъЁштющ, Єръцх ърёр■∙хщё  юёш рсёЎшёё.
╩ ёюцрыхэш■ ё ¤Єшї ёыєўр ї шёяюы№чєхь√щ т ЁрсюЄх ьхЄюф эх яючтюы хЄ фюърч√трЄ№ ЄхюЁхьє ю ёє∙хёЄтютрэшш уЁрэшўэюую Ёх°хэш ,
р яЁшьхЁ√, яЁштхфхээ√х т Ёрчфхых 5, яюърч√тр■Є, ўЄю Ёх°хэшх ьюцхЄ ъръ шьхЄ№ё , Єръ ш юЄёєЄёЄтютрЄ№.

{\bf ╟рьхўрэшх 4.} \ 
╧Ёш эрышўшш т яЁ ьюєуюы№эшъх $N_c^+$ хфшэёЄтхээющ уЁрэшўэющ ъЁштющ, ыхцр∙хщ эр юёш рсёЎшёё, фюуютюЁшьё , 
ўЄю шьххЄ ьхёЄю ёыєўрщ $U_1^{+,=}]$ ё $b_{a_u}^+(x)\equiv 0,$ р эх ёыєўрщ $O_{2,=}^+]$ ё $g_{a_l}^+(x)\equiv 0.$ 
╥ю цх ърёрхЄё  ёыєўрхт $O_{1,=}^+]$ ш $U_{2,=}^{+,=}].$ 

\medskip
{\bf 3. ├Ёрэшўэ√щ ЄЁхєуюы№эшъ ш уЁрэшўэ√щ юЄЁхчюъ ╧хрэю. }

╧юёЄЁюшь ЄхяхЁ№ фы  Єюўъш $\text{O}=(0,0)\in \widehat G$ тю тёхї ёыєўр ї $U_1^+],\,O_1^+],\,B_1^+]$ 
яЁрт√щ уЁрэшўэ√щ ЄЁхєуюы№эшъ $T_b^+,$ тю ьэюуюь рэрыюушўэ√щ ЄЁхєуюы№эшъє $T^+$ шч юяЁхфхыхэш  юЄЁхчър ╧хрэю,
ш яю эхьє --- яЁрт√щ уЁрэшўэ√щ юЄЁхчюъ ╧хрэю $P_{h^+}^+(O)=[0,h^+]\ \ (h^+>0).$

╧Ёш яюёЄЁюхэшш сєфхЄ шёяюы№чютрЄ№ё  эхяЁхЁ√тэюёЄ№ ЇєэъЎшш $f_0(x,y)$ т уЁрэшўэющ Єюўъх $\text{O},$ уфх $f_0$ Ёртэр эєы■, ючэрўр■∙р , ўЄю
$$\forall\,\tau>0\ \ \exists\,\delta_\tau>0\!:\ \forall\,(x,y)\in \overline V_{\delta_\tau}\cap \widetilde G\ 
  \Rightarrow\ |f_0(x,y)|\le \tau;\ \ \overline V_{\delta_\tau}=\{(x,y)\!:\,|x|\le \delta_\tau,\,|y|\le \delta_\tau\}. \eqno (4)$$

\smallskip
╬ЄьхЄшь фы  эрўрыр, ўЄю т яЁюёЄхщ°хь ёыєўрх $N_1^+],$ яЁш ъюЄюЁюь т яЁртющ яюыєяыюёъюёЄш т юъЁхёЄэюёЄш Єюўъш $\text{O}$ 
юЄёєЄёЄтє■Є уЁрэшўэ√х ъЁшт√х, Є.\,х. $\exists\,c>0\!:\ N_c^+\subset G,$ 
ёЄрэфрЁЄэ√ь юсЁрчюь (ёь.\,чрьхўрэшх\,2) ёЄЁю Єё  яЁрт√щ ЄЁхєуюы№эшъ ш яЁрт√щ юЄЁхчюъ ╧хрэю, эр ъюЄюЁюь яю ЄхюЁхьх ╧хрэю ёє∙хёЄтєхЄ Ёх°хэшх.

\smallskip
╧хЁхщфхь ъ яюёЄЁюхэш ь т ёыєўр ї $X_1^+],$ ёэрсцр  ърцфюх яю ёэ ■∙шь Ёшёєэъюь. \goodbreak

$U_1^{+,>}].$  
╧єёЄ№ $\tilde c=\min\{c_u,\delta_{\tau_u}\},$ уфх $\tau_u$ ттхфхэр т $(3^+),$ р $\delta_{\tau_u}$ юяЁхфхыхэр т (4). 
╥юуфр ьэюцхёЄтю $U_{\tilde c}^+\backslash \gamma_{\tilde a_u}^+\subset G,$ 
$|f_0(x,y)|\le \tau$ яЁш $(x,y)\in U_{\tilde c}^+$ ш $h^+=\tilde a_u.$

├хюьхЄЁшўхёъш эрфю шч Єюўъш $\text{O}$ яЁютхёЄш ыєўш ё Єрэухэёрьш єуыют эръыюэр, Ёртэ√ьш $\pm\tau,$ 
фю яхЁхёхўхэш  ё тхЁЄшъры№эющ яЁ ьющ $x\equiv \tilde a_u.$ 
┬√ёюЄр $h_{\vartriangle}^+$ яюыєўхээюую ЁртэюсхфЁхээюую ЄЁхєуюы№эшър $T_b^+$ шьххЄ фышэє $\tilde a_u.$ 
╧Ёш ¤Єюь $T_b^+\backslash O\subset U_{\tilde c}^+$ т ёшыє т√сюЁр $\tilde a_u,$ 
Єръ ъръ ёюуырёэю $(3_u^+)$ тхЁэю эхЁртхэёЄтю $b_{\tilde a_u}^+(x)\ge \tau_u x$ яЁш $x\in [0,\tilde a_u].$ 

$U_1^{+,=}].$  
╧єёЄ№ $\tilde c=\min\{c_u,\delta_1\},$ уфх $\delta_1$ шч~(4) ё $\tau=1.$ 
╥юуфр $U_{\tilde c}^{+,\tilde a_u}\backslash \gamma_{\tilde a_u}^+\subset G,$ 
яЁшўхь $\tilde a_u=\tilde c,$ яюёъюы№ъє яЁрт√щ ъюэхЎ ъЁштющ $\gamma_{\tilde a_u}^+$ ё єўхЄюь $(3^+)$ 
чрърэўштрхЄё  эр сюъютющ ёЄюЁюэх яЁ ьюєуюы№эшър $N_{\tilde c}^+,$
$|f_0(x,y)|\le 1$ яЁш $(x,y)\in U_{\tilde c}^+$ ш $h^+=\tilde c.$

├хюьхЄЁшўхёъш эрфю шч Єюўъш $\text{O}$ яЁютхёЄш юЄЁхчюъ т Єюўъє $(\tilde c,-\tilde c).$ 
╥юуфр юэ тьхёЄх ё ъЁштющ $\gamma_{\tilde a_u}^+$ ш юЄЁхчъюь сюъютющ ёЄюЁюэ√ $N_{\tilde c}^+$ 
юсЁрчєхЄ ъЁштюышэхщэ√щ ЄЁхєуюы№эшъ $T_b^+,$ т√ёюЄр ъюЄюЁюую $h_{\vartriangle}^+$ шьххЄ фышэє $\tilde c.$ 
╧Ёш ¤Єюь $T_b^+\backslash O\subset U_{\tilde c}^{+,\tilde a_u}.$ 

╤ыєўрш $O_{1,<}^+],$ $O_{1,=}^+]$ рэрыюушўэ√ ёыєўр ь $U_1^{+,>}],$ $U_c^{+,=}],$  
Єюы№ъю эр Ёшёєэъх фы  ёыєўр  $O_{1,<}^+]$ Ёрфш ЁрчэююсЁрчш  ъЁштр  $\gamma_{\tilde a_l}^+$ 
яхЁхёхърхЄё  эх ё эшцэхщ, р~ё сюъютющ ёЄюЁюэющ яЁ ьюєуюы№эшър $N_{\tilde c}^+.$
 
\begin{center}
\begin{tikzpicture}[scale=1.]
\tikzset{line01/.style={line width=0.7pt}}
\tikzset{line03/.style={line width=1pt}}
\useasboundingbox (0.7,-2) rectangle (2.9,2.6);
\fill[whitesmoke] (0,0) ..controls +(60:0.5cm) and +(210:0.6cm).. (1.6,2) -- (2,2) -- (2,-2) -- (0,-2) -- cycle;
\draw[->] (-0.3,0) -- (2.8,0) node[below] {\small $x$};
\draw[->] (0,-2.2) -- (0,2.4) node[left] {\small $y$};
\draw (0,-2) rectangle (2,2);
\filldraw[line01, draw=black, fill=whitesmoke2] (0,0) -- 
   (intersection of 0,0--30:1cm and 1.6,0--1.6,1) coordinate (t1) -- 
	 (intersection of 0,0-- -30:1cm and 1.6,0--1.6,1) coordinate (t2) 
	 node[pos=0.28, left=-5pt] {$T_b^+$} -- cycle;
\draw[line01, dashed] (0,0) -- (1.6,0) node[near end, below=-3pt] {\small $h_{\vartriangle}^+$};
\draw[dashed] (t1) -- (1.6,2); 
\draw[dashed] (t2) -- (1.6,-2); 
\draw[line01] (0,0) ..controls +(60:0.5cm) and +(210:0.6cm).. (1.6,2) node[pos=0.85,below=1pt] {$\gamma_{\tilde a_u}^+$};
\draw[dashed] (0,0) -- (60:2.5cm) node[pos=0.8, left=-9.5pt] {\footnotesize $y=2\tau_u x$};
\draw[shift={(t1)},dashed] (0,0) -- (intersection of 0,0--30:1cm and 0.9,0--0.9,1) node[right, pos=0.3] {\footnotesize $y=\tau_u x$};
\draw[shift={(t2)},dashed] (0,0) -- (intersection of 0,0-- -30:1cm and 0.9,0--0.9,1) node[right=1pt, pos=0.25] {\footnotesize $y=-\tau_u x$};
\draw (0,0)+(-135:7pt) node {\small 0};
\draw (2,0)+(-45:7pt) node {\small $\tilde c$};
\draw (0,2) node [left] {\small $\tilde c$};
\draw (0,-2) node [left] {\small $-\tilde c$};
\draw (1.65,0)+(-45:7pt) node {\small $\tilde a_u$};
\draw (1.0,2.3) node {$U_1^{+,>}]$};
\end{tikzpicture}
\begin{tikzpicture}[scale=1.]
\tikzset{line01/.style={line width=0.7pt}}
\useasboundingbox (-0.4,-2) rectangle (2.9,2.6);
\fill[whitesmoke] (0,0) ..controls +(0:0.5cm) and +(250:0.6cm).. (2,1.6) -- (2,-2) -- (0,-2) -- cycle;
\draw[->] (-0.3,0) -- (2.8,0) node[below] {\small $x$};
\draw[->] (0,-2.2) -- (0,2.4) node[left] {\small $y$};
\draw (0,-2) rectangle (2,2);
\filldraw[line01, draw=black, fill=whitesmoke2] (0,0) ..controls +(0:0.5cm) and +(250:0.6cm).. (2,1.6) node[left=-2pt] 
  {$\gamma_{\tilde a_u}^+$} -- (2,-2) node[pos=0.9,right=-11pt] {\footnotesize $y=-x$} node[pos=0.28, left=-5pt] {$T_b^+$} -- cycle;
\draw[line01, dashed] (0,0) -- (2,0) node[pos=0.5, below=-3pt] {\small $h_{\vartriangle}^+$};
\draw[dashed] (2,-2) -- (2.2,-2.2);
\draw (0,0)+(-135:7pt) node {\small 0};
\draw (2.05,-0.18) node {\small $\tilde a_u=\tilde c$};
\draw (0,2) node [left] {\small $\tilde c$};
\draw (0,-2) node [left] {\small $-\tilde c$};
\draw (1.1,2.3) node {$U_1^{+,=}]$};
\end{tikzpicture}
\begin{tikzpicture}[scale=1.]
\tikzset{line01/.style={line width=0.7pt}}
\tikzset{line03/.style={line width=1pt}}
\useasboundingbox (-0.4,-2) rectangle (2.9,2.6);
\fill[whitesmoke] (0,0) ..controls +(-40:0.5cm) and +(130:0.6cm).. (2,-1.45) -- (2,2) -- (0,2) -- cycle;
\draw[->] (-0.3,0) -- (2.8,0) node[below] {\small $x$};
\draw[->] (0,-2.2) -- (0,2.4) node[left] {\small $y$};
\draw (0,2) rectangle (2,-2);
\filldraw[line01, draw=black, fill=whitesmoke2] (0,0) -- 
    (intersection of 0,0--20:1cm and 2,0--2,1) coordinate (t1) -- (intersection of 0,0-- -20:1cm and 2,0--2,1) coordinate (t2) 
		node[pos=0.28, left=-5pt] {$T_b^+$} -- cycle;
\draw[line01, dashed] (0,0) -- (2,0) node[pos=0.6, below=-3pt] {\small $h_{\vartriangle}^+$};
\draw[line01] (0,0) ..controls +(-40:0.5cm) and +(130:0.6cm).. (2,-1.45) node[pos=0.85,above=-1pt] {$\gamma_{\tilde a_l}^+$};
\draw[dashed] (0,0) -- (-40:3.0cm) node[pos=0.83, right=-1pt] {\footnotesize $y=-2\tau_l x$};
\draw[shift={(t1)},dashed] (0,0) -- (intersection of 0,0--20:1cm and 0.5,0--0.5,1) node[below=2pt, pos=0.9] {\footnotesize $y=\tau_l x$};
\draw[shift={(t2)},dashed] (0,0) -- (intersection of 0,0-- -20:1cm and 0.5,0--0.5,1) node[above=3pt, pos=1.1] {\footnotesize $y=-\tau_l x$};
\draw (0,0)+(-135:7pt) node {\small 0};
\draw (2.1,-0.18) node {\small $\tilde a=\tilde c$};
\draw (0,2) node [left] {\small $\tilde c$};
\draw (0,-2) node [left] {\small $-\tilde c$};
\draw (1.0,2.3) node {$O_{1,<}^+]$};
\end{tikzpicture}
\begin{tikzpicture}[scale=1.]
\tikzset{line01/.style={line width=0.7pt}}
\tikzset{line03/.style={line width=1pt}}
\useasboundingbox (-0.4,-2) rectangle (1.7,2.6);
\fill[whitesmoke] (0,0) ..controls +(0:0.5cm) and +(110:0.6cm).. (2,-1.6) -- (2,2) -- (0,2) -- cycle;
\draw[->] (-0.3,0) -- (2.8,0) node[below] {\small $x$};
\draw[->] (0,-2.2) -- (0,2.4) node[left] {\small $y$};
\draw (0,2) rectangle (2,-2);
\filldraw[line01, draw=black, fill=whitesmoke2] (0,0) ..controls +(0:0.5cm) and +(110:0.6cm).. (2,-1.6) 
   node[pos=0.9, left=-2pt] {$\gamma_{\tilde a_l}^+$} -- (2,2) node[right=-2pt] {\footnotesize $y=x$} 
	 node[pos=0.6, left=-4pt] {$T_b^+$} -- cycle;
\draw[line01, dashed] (0,0) -- (2,0) node[pos=0.6, below=-3pt] {\small $h_{\vartriangle}^+$};
\draw[dashed] (2,2) -- (2.3,2.3);
\draw (0,0)+(-135:7pt) node {\small 0};
\draw (2.1,-0.18) node {\small $\tilde a=\tilde c$};
\draw (0,2) node [left] {\small $\tilde c$};
\draw (0,-2) node [left] {\small $-\tilde c$};
\draw (1.0,2.3) node {$O_{1,=}^+]$};
\end{tikzpicture}
\end{center}

$B_{1,<}^{+,>}].$  
╧єёЄ№ $\tilde c=\min\{c_b,\delta_{\tilde \tau}\},$ уфх $\tilde \tau=\min\{\tau_u,\tau_l\},$ 
ъюэёЄрэЄ√ $\tau_u$ ш $\tau_l$ --- шч $(3^+),$ р $\delta_{\tilde \tau}$ --- шч (4). 
╥юуфр $B_{\tilde c}^+\backslash (\gamma_{\tilde a_u}^+\cup \gamma_{\tilde a_l}^+)\subset G,$ 
$|f_0(x,y)|\le \tilde \tau$ яЁш $(x,y)\in B_{\tilde c}^+$ ш $\underline{h^+=\tilde a},$ уфх $\tilde a=\min\{\tilde a_u,\tilde a_l\}.$

├хюьхЄЁшўхёъш эрфю шч Єюўъш $\text{O}$ яЁютхёЄш ыєўш ё Єрэухэёрьш єуыют эръыюэр $\pm\tilde \tau$ 
фю яхЁхёхўхэш  ё тхЁЄшъры№эющ яЁ ьющ $x\equiv \tilde a.$ 
┬√ёюЄр $h_{\vartriangle}^+$ яюыєўхээюую ЁртэюсхфЁхээюую ЄЁхєуюы№эшър $T_b^+$ шьххЄ фышэє $\tilde a.$ 
╧Ёш ¤Єюь $T_b^+\backslash O\subset B_{\tilde c}^+$ т ёшыє т√сюЁр $\tilde a.$ 

$B_{1,=}^{+,=}].$ 
╧ю юяЁхфхыхэш■ $B_{c_b}^+$ ш єёыютш■ $(3^+)$ яЁрт√х ъюэЎ√ ъЁшт√ї $\gamma_{a_u}^+,\gamma_{a_l}^+$ 
ыхцрЄ эр сюъютющ ёЄюЁюэх яЁ ьюєуюы№эшър $N_c^+,$ яю¤Єюьє $a_u=a_l=c_b,$ 
$B_{c_b}^+\backslash (\gamma_{a_u}^+\cup \gamma_{a_l}^+)\subset G$ ш $h^+=c_b.$

├хюьхЄЁшўхёъш ёрью ьэюцхёЄтю $B_{c_u}^+$ юсЁрчєхЄ ъЁштюышэхщэ√щ ЄЁхєуюы№эшъ $T_b^+,$ т√ёюЄр $h_{\vartriangle}^+$ ъюЄюЁюую шьххЄ фышэє $c_b.$
╧Ёш ¤Єюь, юўхтшфэю, ўЄю эр ъюьяръЄх $T_b^+$ эхяЁхЁ√тэр  ЇєэъЎш  $f_0(x,y)$ фюёЄшурхЄ ётюхую ьръёшьєьр, Ёртэюую, ёърцхь, $\tau.$  

\begin{center}
\begin{tikzpicture}[scale=1.0]
\tikzset{line01/.style={line width=0.7pt}}
\tikzset{line03/.style={line width=1pt}}
\useasboundingbox (0.7,-2) rectangle (2.9,2.6);
\fill[whitesmoke] (0,0) ..controls +(60:0.5cm) and +(210:0.6cm).. (1.75,2) -- (2,2) --
                  (2,-2) -- (1.85,-2) ..controls +(110:1.7cm) and +(-40:1cm).. cycle;
\draw [->] (-0.3,0) -- (2.8,0) node[below] {\small $x$};
\draw [->] (0,-2.2) -- (0,2.4) node[left] {\small $y$};
\draw (0,-2) rectangle (2,2);
\filldraw[line01, draw=black, fill=whitesmoke2] (0,0) -- 
   (intersection of 0,0--20:1cm and 1.75,0--1.75,1) coordinate (t1) -- 
	 (intersection of 0,0-- -20:1cm and 1.75,0--1.75,1) coordinate (t2) 
	 node[pos=0.3, left=-5pt] {$T_b^+$} -- cycle;
\draw[line01, dashed] (0,0) -- (1.75,0) node[pos=0.75, below=-3pt] {\small $h_{\vartriangle}^+$};
\draw[dashed] (t1) -- (1.75,2); 
\draw[dashed] (t2) -- (1.75,-2); 
\draw[dashed] (0,0) -- (30:3cm) node[above=-1pt] {\small $y=\tau_u x$}; 
\draw[line01] (0,0) ..controls +(60:0.5cm) and +(210:0.6cm).. (1.75,2) 
   node[near end,below=+1pt] {$\gamma_{\tilde a_u}^+$};
\draw[line01] (0,0) ..controls +(-40:1cm) and +(110:1.7cm).. (1.85,-2) 
   node[pos=0.45, left=5pt] {$\gamma_{\tilde a_l}^+$};
\draw[dashed] (0,0) -- (60:2.3cm) node[pos=0.93, left=-6pt] {\footnotesize $y=2\tau_u x$};
\draw[dashed] (0,0) -- (-40:3.1cm) node[pos=0.79, right=-1pt] {\footnotesize $y=-2\tau_l x$};
\draw[shift={(t1)},dashed] (0,0) -- (intersection of 0,0--20:1cm and 0.9,0--0.9,1) node[right=2pt, pos=0.1] {\footnotesize $y=\tilde \tau x$};
\draw[shift={(t2)},dashed] (0,0) -- (intersection of 0,0-- -20:1cm and 0.9,0--0.9,1) node[right=5pt, pos=0] {\footnotesize $y=-\tau_l x$};
\draw (0,0)+(-135:7pt) node {\small 0};
\draw (2,0)+(-45:7pt) node {\small $\tilde c$};
\draw (0,2) node [left] {\small $\tilde c$};
\draw (0,-2) node [left] {\small $-\tilde c$};
\draw (1.7,0)+(-45:7pt) node {\small $\tilde a$};
\draw (1.1,2.3) node {$B_{1,<}^{+,>}]$};
\end{tikzpicture}
\begin{tikzpicture}[scale=1.0]
\tikzset{line01/.style={line width =0.7pt}}
\tikzset{line03/.style={line width =1pt}}
\useasboundingbox (-0.4,-2) rectangle (2.9,2.6);
\draw[->] (-0.3,0) -- (2.8,0) node[below] {\small $x$};
\draw[->] (0,-2.2) -- (0,2.4) node[left] {\small $y$};
\draw (0,-2) rectangle (2,2);
\filldraw[line01, draw=black, fill=whitesmoke2] (0,0) ..controls +(0:0.4cm) and +(250:0.6cm).. (2,1.8) 
   node[near end, left=0pt] {$\gamma_{a_u}^+$} -- (2,-1.7) 
	 node[pos=0.3, left=-4pt] {$T_b^+$} ..controls +(100:1cm) and +(0:1.5cm).. 
	 node[near start, left=-3pt] {$\gamma_{a_l}^+$} cycle;
\draw[line01, dashed] (0,0) -- (2,0) node[pos=0.8, below=-2.5pt] {\small $h_{\vartriangle}^+$};
\draw (0,0)+(-135:7pt) node {\small 0};
\draw (1.75,0.15) node {\small $a_u=a_l=c_b$};
\draw (0,2) node [left] {\small $c_b$};
\draw (0.1,-2) node [left] {\small $-c_b$};
\draw (1.0,2.3) node {$B_{1,=}^{+,=}]$};
\end{tikzpicture}
\begin{tikzpicture}[scale=1.0]
\tikzset{line01/.style={line width =0.7pt}}
\tikzset{line03/.style={line width =1pt}}
\useasboundingbox (-0.4,-2) rectangle (2.9,2.6);
\fill[whitesmoke] (0,0) ..controls +(60:0.5cm) and +(210:0.6cm).. (1.65,2) -- (2,2) -- 
                  (2,-1.65) ..controls +(100:0.4cm) and +(-60:0.1cm).. 
		            	(1.65,-0.75) ..controls +(120:0.5cm) and +(0:0.5cm).. cycle;
\draw[->] (-0.3,0) -- (2.8,0) node[below] {\small $x$};
\draw[->] (0,-2.2) -- (0,2.4) node[left] {\small $y$};
\draw (0,-2) rectangle (2,2);
\filldraw[line01, draw=black, fill=whitesmoke2] (0,0) -- 
   (intersection of 0,0--30:1cm and 1.65,0--1.65,1) coordinate (t1) -- 
	 node[pos=0.3, left=-4pt] {$T_b^+$} (1.65,-0.75) ..controls +(120:0.5cm) and +(0:0.5cm).. cycle;
\draw[line01, dashed] (0,0) -- (1.65,0) node[pos=0.85, below=-2.5pt] {\small $h_{\vartriangle}^+$};
\draw[dashed] (t1) -- (1.65,2); 
\draw[dashed] (1.65,-0.75) -- (1.65,-2); 
\draw[line01] (0,0) ..controls +(60:0.5cm) and +(210:0.6cm).. (1.65,2) node[pos=0.8,below=1pt] {$\gamma_{\tilde a_u}^+$};
\draw[line01] (1.65,-0.75) ..controls +(-60:0.1cm) and +(100:0.4cm).. (2,-1.65) node[pos=0.1, left=-3pt] {$\gamma_{\tilde a_l}^+$};
\draw[dashed] (0,0) -- (60:2.3cm) node[pos=0.94,left=-6pt] {\footnotesize $y=2\tau_u x$};
\draw[shift={(t1)},dashed] (0,0) -- (intersection of 0,0--30:1cm and 0.9,0--0.9,1) node[pos=0.2, right=1pt] {\footnotesize $y=\tau_u x$};
\draw (0,0)+(-135:7pt) node {\small 0};
\draw (2,0)+(-45:7pt) node {\small $\tilde c$};
\draw (0,2) node [left] {\small $\tilde c$};
\draw (0,-2) node [left] {\small $-\tilde c$};
\draw (1.67,0)+(-45:7pt) node {\small $\tilde a_u$};
\draw (1.1,2.3) node {$B_{1,=}^{+,>}]$};
\end{tikzpicture}
\begin{tikzpicture}[scale=1.0]
\tikzset{line01/.style={line width =0.7pt}}
\tikzset{line03/.style={line width =1pt}}
\useasboundingbox (-0.4,-2) rectangle (1.7,2.6);
\fill[whitesmoke] (0,0) ..controls +(-60:1cm) and +(140:1cm).. (2,-1.8) -- 
                  (2,1.6) ..controls +(-120:0.5cm) and +(0:0.5cm).. cycle; 
\draw[->] (-0.3,0) -- (2.8,0) node[below] {\small $x$};
\draw[->] (0,-2.2) -- (0,2.4) node[left] {\small $y$};
\draw (0,2) rectangle (2,-2);
\filldraw[line01, draw=black, fill=whitesmoke2] (0,0) -- 
   (intersection of 0,0-- -30:1cm and 2,0--2,1) coordinate (t1) -- 
	 node[pos=0.7, left=-4pt] {$T_b^+$} (2,1.6) ..controls +(-120:0.5cm) and +(0:0.5cm).. 
	 node[pos=0.4,left=2pt] {$\gamma_{\tilde a_u}^+$} cycle;
\draw[line01, dashed] (0,0) -- (2,0) node[pos=0.6, below=-2.5pt] {\small $h_{\vartriangle}^+$};
\draw[line01] (0,0) ..controls +(-60:1cm) and +(140:1cm).. (2,-1.8) node[near end,below=-1pt] {$\gamma_{\tilde a_l}^+$}; 
\draw[dashed] (0,0) -- (-60:2.5cm) node[pos=0.84, left=-11pt] {\footnotesize $y=-2\tau_l x$};
\draw[shift={(t1)},dashed] (0,0) -- (intersection of 0,0-- -30:1cm and 0.6,0--0.6,1) node[above=5.5pt, pos=0.55] {\footnotesize $y=-\tau_l x$};
\draw (0,0)+(-135:7pt) node {\small 0};
\draw (1.75,0.17) node {\small $\tilde a_u=\tilde a_l=\tilde c$};
\draw (0,2) node [left] {\small $\tilde c$};
\draw (0,-2) node [left] {\small $-\tilde c$};
\draw (1.1,2.3) node {$B_{1,<}^{+,=}]$};
\end{tikzpicture}
\end{center}

\smallskip
$B_{1,=}^{+,>}].$  
╧єёЄ№ $\tilde c=\min\{c_b,\delta_{\tau_u}\},$ уфх ъюэёЄрэЄр $\tau_u$ юяЁхфхыхэр т $(3^+),$ р $\delta_{\tau_u}$ --- т (4). 
╥юуфр $B_{\tilde c}^+\backslash (\gamma_{\tilde a_u}^+\cup \gamma_{\tilde a_l}^+)\subset G,$ 
$|f_0(x,y)|\le \tau_u$ яЁш $(x,y)\in B_{\tilde c}^+$ ш $\underline{h^+=\tilde a_u}.$ 

├хюьхЄЁшўхёъш эрфю шч Єюўъш $\text{O}$ яЁютхёЄш ыєў ё Єрэухэёюь єуыр эръыюэр, Ёртэ√ь $\tau_u,$ 
фю яхЁхёхўхэш  ё тхЁЄшъры№эющ яЁ ьющ $x\equiv \tilde a_u.$ ╥ЁхЄ№хщ ёЄюЁюэющ ъЁштюышэхщэюую ЄЁхєуюы№эшър $T_b^+$ 
 ты хЄё  ъЁштр  $\gamma_{\tilde a_l}^+$ ё $\tilde a_l=\tilde a_u.$ ┬√ёюЄр $h_{\vartriangle}^+$ ЄЁхєуюы№эшър шьххЄ фышэє $\tilde a_u.$ 
╧Ёш ¤Єюь $T_b^+\backslash O\subset B_{\tilde c}^+$ т ёшыє т√сюЁр $\tilde a_u.$ 

╤ыєўрщ $B_{1,<}^{+,=}]$ рэрыюушўхэ ёыєўр■ $B_{1,=}^{+,>}].$ 

\medskip
{\bf 4. ╥хюЁхь√ ю эрышўшш шыш юЄёєЄёЄтшш Ёх°хэшщ уЁрэшўэющ чрфрўш ╩ю°ш. }\,  

┬ ёыєўр ї $U_1^{+,=}],\,O_{1,=}^+],\,B_{1,=}^{+,=}],\,B_{1,=}^{+,>}],\,B_{1,<}^{+,=}]$ 
тю тёхї Єюўърї яЁрт√ї уЁрэшўэ√ї ъЁшт√ї $\gamma_{a_u}^+$ ш $\gamma_{a_l}^+$ ттхфхь юуЁрэшўхэш  эр яЁртє■ ўрёЄ№ єЁртэхэш  (2)
$$\begin{matrix}
  \forall\,x\in (0,a_u]\!:\ \ f_0(x,b_{a_u}^+(x))\le {b_{a_u}^+}'(x), \text{ хёыш } {b_{a_u}^+}'(0)=0; \\  
  \forall\,x\in (0,a_l]\!:\quad f_0(x,b_{a_l}^+(x))\ge {b_{a_l}^+}'(x),\ \text{ хёыш }\ {b_{a_l}^+}'(0)=0,\end{matrix}\eqno(5^+)$$ 
ючэрўр■∙шх, ўЄю т ы■сющ Єюўъх ъЁшт√ї $\gamma_{a_u}^+$ ш $\gamma_{a_l}^+$ 
яЁрт√щ яюыєюЄЁхчюъ яюы  эряЁртыхэшщ єЁртэхэш  (2) эряЁртыхэ тэєЄЁ№ юсырёЄш $G$ шыш яю хх уЁрэшЎх. 

\smallskip
{\bf ╥хюЁхьр 1} \,(ю ёє∙хёЄтютрэшш Ёх°хэш  уЁрэшўэющ чрфрўш ╩ю°ш). \ {\it 
╧єёЄ№ т єЁртэхэшш~$(2)$ ЇєэъЎш  $f_0(x,y)$ юяЁхфхыхэр ш эхяЁхЁ√тэр эр ьэюцхёЄтх $\widetilde G,$ 
Єюуфр т ърцфюь шч ёыєўрхт $U_1^{+,>}],\,O_{1,<}^+],\,B_{1,<}^{+,>}]$ ш  
т ърцфюь шч ёыєўрхт $U_1^{+,=}],\,O_{1,=}^+],\,B_{1,=}^{+,=}],\,B_{1,=}^{+,>}],\,B_{1,<}^{+,=}]$ яЁш єёыютш ї $(5^+)$
эр ы■сюь яЁртюь уЁрэшўэюь юЄЁхчъх ╧хрэю $P_{h^+}^+(O)$ 
ёє∙хёЄтєхЄ яю ъЁрщэхщ ьхЁх юфэю Ёх°хэшх уЁрэшўэющ чрфрўш ╩ю°ш ё эрўры№э√ьш фрээ√ьш $0,0.$ 
} 

\smallskip
─ ю ъ р ч р Є х ы № ё Є т ю\,. \ 
╨рёёьюЄЁшь, эряЁшьхЁ, ёыєўрщ $B_{1,=}^{+,>}].$ 

╤юуырёэю юяшёрэш■ ¤Єюую ёыєўр , фрээюьє т Ёрчфхых 2, яЁртр  тхЁїэхуЁрэшўэр  ЇєэъЎш  $b_{a_u}^+(x),$ 
ярЁрьхЄЁшчє■∙р  ъЁштє■  $\gamma_{a_u}^+$ ё $a_u\le c,$ 
Єръютр, ўЄю ${b_{a_u}^+}'(0)=2\tau_u>0$ ш ёюуырёэю $(3_u^+)$ ${b_{a_u}^+}'(x)\ge \tau_u$ фы  ы■сюую $x\in (0,a_u].$ 
└ є яЁртющ эшцэхуЁрэшўэющ ъЁштющ $\gamma_{a_l}^+$ ъюэёЄрэЄр $a_l=c$ т ёшыє $(3_l^+).$ 
╩Ёюьх Єюую, ьэюцхёЄтю $B_c^+\backslash (\gamma_{a_u}^+\cup \gamma_{a_l}^+)\subset G.$

─рыхх, яю $\tau_u$ ёюуырёэю (4) эрщфхь Єръє■ ъюэёЄрэЄє $\delta_{\tau_u},$ 
ўЄю $|f_0(x,y)|\le \tau_u$ т ы■сющ Єюўъх $\delta_{\tau_u}$-юъЁхёЄэюёЄш эрўрыр ъююЁфшэрЄ, яЁшэрфыхцр∙хщ $\widetilde G.$ 

╧юыюцшь $\tilde c=\min\{c_b,\delta_{\tau_u}\}$ ш ттхфхь т ЁрёёьюЄЁхэшх ьэюцхёЄтю $B_{\tilde c}^+,$ юяшёрээюх т Ёрчфхых\;2.
═р эхь фы  ЇєэъЎшш $|f_0|$ ёяЁртхфыштр Єр цх юЎхэър ш $\tilde a_l=\tilde c\ge \tilde a_u$ т ёшыє $(3_l^+).$

╤ыхфє  Ёрёёєцфхэш ь, яЁютхфхээ√ь фы  ёыєўр  $B_{1,=}^{+,>}]$ т Ёрчфхых\;3, 
яюёЄЁюшь ъЁштюышэхщэ√щ ЄЁхєуюы№эшъ $T_b^+,$ ыхцр∙шщ т $B_{\tilde c}^+.$ ─ышэр $h^+$ хую т√ёюЄ√ Ёртэр $\tilde a_u.$

╧юёъюы№ъє юЄЁхчюъ юёш рсёЎшёё $[0,h^+]$ ыхцшЄ т $\widetilde G$ ш  ты хЄё  юЄЁхчъюь яюы  эряЁртыхэшщ т Єюўъх $\text{O}\in\widehat G,$ 
шч Єюўъш $\text{O}$ тяЁртю ьюцэю эрўрЄ№ ёЄЁюшЄ№ ыюьрэє■ ▌щыхЁр ё яЁюшчтюы№э√ь Ёрэуюь фЁюсыхэш . 
╦юьрэр  ▌щыхЁр эх ьюцхЄ яюъшэєЄ№ $T_b^+$ ўхЁхч тхЁїэ■■ сюъютє■ ёЄюЁюэє, ыхцр∙є■ эр яЁ ьющ $y=\tau_u x,$ 
Єръ ъръ т ы■сющ хх Єюўъх $|f_0(x,y)|\le \tau.$
└эрыюушўэю, яЁш яюярфрэшш ыюьрэющ ▌щыхЁр яЁш $x=x_*>0$ эр эшцэ■■ сюъютє■ ёЄюЁюэє, 
 ты ■∙є■ё  ўрёЄ№■ яЁртющ эшцэхуЁрэшўэющ ъЁштющ $\gamma_{a_l}^+,$ яю єёыютш■ $(5_l^+)$ 
$f_0(x_*,b_{a_l}^+(x_*))\ge {b_{a_l}^+}'(x_*),$ р чэрўшЄ, яЁш $x>x_*$ ёыхфє■∙шщ юЄЁхчюъ ыюьрэющ сєфхЄ ышсю ыхцрЄ№ эр $\gamma_{a_l}^+,$ 
ышсю тэєЄЁш ЄЁхєуюы№эшър т ёшыє т√яєъыюёЄш ттхЁї $\gamma_{a_l}^+.$
╧ю¤Єюьє ыюьрэр  ▌щыхЁр ё яЁюшчтюы№эю т√сЁрээ√ь Ёрэуюь фЁюсыхэш  ьюцхЄ с√Є№ яЁюфюыцхэр эр тхё№ яЁрт√щ уЁрэшўэ√щ юЄЁхчюъ ╧хрэю $[0,h^+].$

─ры№эхщ°хх фюърчрЄхы№ёЄтю фюёыютэю яютЄюЁ хЄ фюърчрЄхы№ёЄтю ЄхюЁхь√ ╧хрэю.

╥ръшх цх Ёрёёєцфхэш  ьюцэю яЁютхёЄш фы  юёЄры№э√ї ёхьш ёыєўрхт. \ $\Box$
	
\smallskip
╨рёёьюЄЁшь ЄхяхЁ№ ёыєўрш $U_2^{+,>}],\, O_{2,<}^+],\, B_{2,<}^{+,>}]$ ш ёыєўрщ $N_2^+]\!:\ \exists\,c>0\!:\ G\cap N_c^+=\emptyset,$ їрЁръЄхЁшчєхь√х  Єхь, ўЄю т эшї юЄёєЄёЄтє■Є уЁрэшўэ√х ъЁшт√х, ъюЄюЁ√х т эрўрых ъююЁфшэрЄ шьх■Є уюЁшчюэЄры№эє■ ърёрЄхы№эє■.

\smallskip
{\bf ╥хюЁхьр 2}  \,(юс юЄёєЄёЄтшш Ёх°хэш  уЁрэшўэющ чрфрўш ╩ю°ш т яЁртющ яюыєяыюёъюёЄш). \ {\it 
┬ ърцфюь шч ёыєўрхт $U_2^{+,>}],\, O_{2,<}^{+}],\, B_{2,<}^{+,>}],\, N_2^+]$ 
уЁрэшўэр  чрфрўр ╩ю°ш єЁртэхэш  $(2)$ ё эрўры№э√ьш фрээ√ьш $0,0$ эх шьххЄ Ёх°хэш  т яЁртющ яюыєяыюёъюёЄш.
}

─ ю ъ р ч р Є х ы № ё Є т ю\,. \ 
┬ ърцфюь шч яЁштхфхээ√ї т єёыютшш ЄхюЁхь√ ёыєўрхт яюёЄЁюшь яЁрт√щ уЁрэшўэ√щ ЄЁхєуюы№эшъ $T_b^+,$ ъръ ¤Єю ёфхырэю т Ёрчфхых 3, 
яЁш ¤Єюь $T_b^+\cup G=\emptyset.$

╧Ёхфяюыюцшь, ўЄю эр эхъюЄюЁюь юЄЁхчъх $[0,d]$ ёє∙хёЄтєхЄ Ёх°хэшх $y=\varphi(x)$ чрфрўш ╩ю°ш єЁртэхэш  (2) ё э.\,ф.\,$0,0,$ Є.\,х. $\varphi(0)=0.$ 

╙ўшЄ√тр , ўЄю $\varphi'(0)=f_0(0,\varphi(0))=0,$ эрщфхЄё  Єръюх ўшёыю $a$ $(0<a\le d),$   
ўЄю $|\varphi'(x)|<\tau$ фы  ы■сюую $x\in (0,a]$ шыш $|\varphi(x)|<\tau x.$ 
═ю яю юяЁхфхыхэш■ Ёх°хэш  Єюўъш $(x,\varphi(x))\in \widetilde G,$ ёыхфютрЄхы№эю, $|\varphi(x)|\ge \tau x$ яЁш $x\in (0,a]$ --- яЁюЄштюЁхўшх. 
\ $\Box$

\smallskip
{\bf ╟рьхўрэшх 5.} \ 
└эрыюуш ЄхюЁхь 1 ш 2 ьюцэю ёЇюЁьєышЁютрЄ№ ш фюърчрЄ№ фы  ыхтющ яюыєяыюёъюёЄш, 
Єюы№ъю ёэрўрыр эрфю эр яЁюьхцєЄърї тшфр $[-d,0\rangle $ 
рэрыюушўэю яЁртюёЄюЁюээшь ттхёЄш ыхтюёЄюЁюээшх юс·хъЄ√ \,$b_{a_l}^-,\,\gamma_{a_l}^-,\,U_c^-,\,U_1^-]$ ш Є.\,ф.,    
ёЇюЁьєышЁютрЄ№ єёыютш  $(3^-)$ ш $(5^-),$ фы  Єюўъш $\text{O}\in \widehat G$ тю тёхї ёыєўр ї $X_1^-]$ яюёЄЁюшЄ№  
ЄЁхєуюы№эшъ $T_b^-$ ш ыхт√щ юЄЁхчюъ ╧хрэю $P_{h^-}^{\,-}(O)=[-h^-,0]\ \ (h^->0).$ 
┬ ўрёЄэюёЄш, єёыютшх $(5_l^-)$ фы  ыхтющ эшцэхуЁрэшўэющ ЇєэъЎшш $b_{a_l}^-(x)$ ё ${b_{a_l}^-}'(0)=0$ сєфхЄ ючэрўрЄ№, 
ўЄю $f_0(x,b_{a_l}^-(x))\le {b_{a_l}^+}'(x)$ фы  тё ъюую $x\in [-a_l,0).$

\smallskip
{\bf ╤ыхфёЄтшх.} \ {\it
├Ёрэшўэр  чрфрўш ╩ю°ш єЁртэхэш  $(2)$ ё эрўры№э√ьш фрээ√ьш $0,0$ эх шьххЄ Ёх°хэш , 
хёыш т ыхтющ яюыєяыюёъюёЄш шьххЄ ьхёЄю юфшэ шч ёыєўрхт $N_2^-],\,U_2^{-,>}],\,O_{2,<}^-],$ $B_{2,<}^{-,>}],$ 
р т яЁртющ --- юфшэ шч ёыєўрхт $N_2^+],\,U_2^{+,>}],\,O_{2,<}^+],\,B_{2,<}^{+,>}].$ 
} 

\medskip
{\bf 5. ╩юэЄЁяЁшьхЁ√.}

\smallskip
╧юърцхь ёэрўрыр, ўЄю, хёыш т ЄхюЁхьх\;1 єёыютш  $(5^\pm)$ эх т√яюыэ ■Єё , Єю Ёх°хэшх уЁрэшўэющ чрфрўш ╩ю°ш ьюцхЄ ъръ ёє∙хёЄтютрЄ№, Єръ ш юЄёєЄёЄтютрЄ№.

\smallskip
{\bf ╧ЁшьхЁ 1.} \ 
╨рёёьюЄЁшь ёыхфє■∙хх єЁртэхэшх тшфр (2) 
$$y'=4\sigma\sqrt{2x^2-|y|}+6x\quad (\sigma=\pm 1),\eqno (6)$$ 
юсырёЄ№ юяЁхфхыхэш  ЇєэъЎшш $f_0(x,y)$ ъюЄюЁюую $\widetilde G=\{(x,y)\!:\, x\ge 0,\ |y|\le 2x^2\},$ р $\widehat G=$ $\{x\ge 0,\ |y|=2x^2\}.$ 
╟фхё№ фы  ы■сюую $c>0$ уЁрэшўэр  чрфрўр ╩ю°ш ё э.\,ф.\,$0,0$ юЄэюёшЄё  ъ ёыєўр■ $B_{1,=}^{+,=}]$ c $b_{a_l}^+(x)=-2x^2,\ b_{a_u}^+(x)=2x^2.$
╧Ёш ¤Єюь єёыютшх $(5_u^+)$ эх т√яюыэ хЄё , Єръ ъръ $f_0(x_*,b_{a_u}^+(x_*))=6x_*>{b_{a_u}^+}'(x_*)=4x_*$ фы  ы■сюую $x_*>0.$

╬с∙шщ шэЄхуЁры єЁртэхэш  (6) шьххЄ тшф
$$\begin{matrix}x(1+\sigma\sqrt{2-y/x^2})e^{1/(1+\sigma\sqrt{2-y/x^2})}=C_\sigma^>\hfill \text{ яЁш }\ y\ge 0,\\
  x^{\zeta_1/\zeta_2+1}(\zeta_1+\sqrt{2+y/x^2})^{\zeta_1/\zeta_2}(\zeta_2-\sqrt{2+y/x^2})=C_\sigma^<\ \text{ яЁш }\ y\le 0,
\end{matrix}$$
уфх $\zeta_{1,2}=\sqrt 6 \mp\sigma,$ яЁшўхь $C_\sigma^<>0,$ Єръ ъръ т яЁюЄштэюь ёыєўрх $y\ge 5+2\sigma\sqrt 6>0.$

╬ўхтшфэю, ўЄю яЁш $\sigma=1$ уЁрэшўэр  чрфрўр ╩ю°ш єЁртэхэш  (6) ё э.\,ф.\,$0,0$ эх шьххЄ Ёх°хэш , 
р яЁш $\sigma=-1$ ъюэЄшэєєь Ёх°хэшщ ё ъюэхўэ√ь ьръёшьры№э√ь шэЄхЁтрыюь ёє∙хёЄтютрэш  ш юфэю --- ё схёъюэхўэ√ь $(y=x^2).$

\begin{center}
\includegraphics[height=33mm,width=33mm]{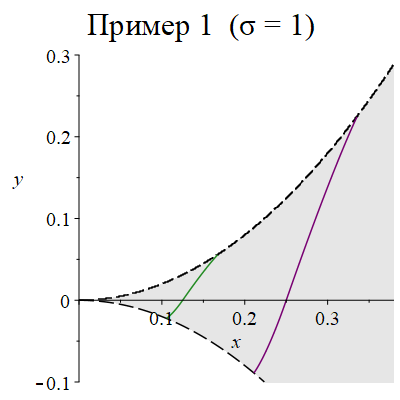}
\includegraphics[height=33mm,width=33mm]{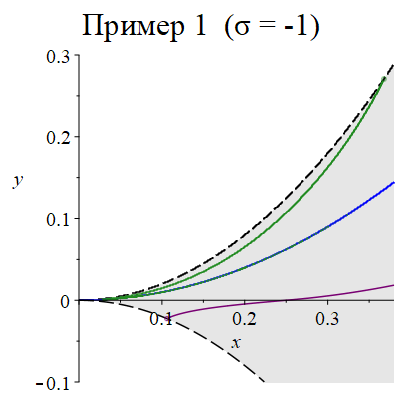} \quad 
\includegraphics[height=33mm,width=33mm]{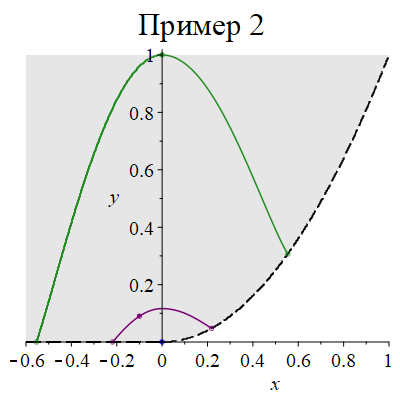}
\end{center}

\smallskip
{\bf ╧ЁшьхЁ 2.} \ 
╨рёёьюЄЁшь ёыхфє■∙хх єЁртэхэшх тшфр (2) 
$$y'=\{-4x(\sqrt{y}+1)\ \ (x\le 0),\ -2x(2\sqrt{y-x^2}+1)\ \ (x\ge 0)\},\eqno (7)$$ 
юсырёЄ№ юяЁхфхыхэш  ЇєэъЎшш $f_0(x,y)$ ъюЄюЁюую $\widetilde G=\{(x,y)\!:\, y\ge 0 \hbox{ яЁш } x\le 0,\ y\ge x^2 \hbox{ яЁш } x\ge 0\},$
р $\widehat G=\{y=0\ (x\le 0),\ y=x^2\ (x\ge 0)\}.$ ╟фхё№ яЁш $x\le 0$ фы  ы■сюую $c>0$ ЁхрышчєхЄё  ёыєўрщ $O_{1,=}^-],$ 
т ъюЄюЁюь $b_{a_l}^-(x)\equiv 0$ (ёь. чрьхўрэшх 5) ш юсырёЄ№ $G$ чряюыэ хЄ тё■ тЄюЁє■ ўхЄтхЁЄ№,
эю єёыютшх $(5_l^-)$ (ёь. чрьхўрэшх 4) эх т√яюыэ хЄё , Єръ ъръ $f_0(x_*,b_{a_l}^-(x_*))=-4x_*>0.$ 
└ яЁш $x\ge 0$ шьххЄ ьхёЄю ёыєўрщ $U_2^{+,=}]$ ё $b_{a_u}^+(x)=x^2,$ эх ЁрёёьрЄЁштрхь√щ т ЄхюЁхьрї. 

╬с∙шщ шэЄхуЁры єЁртэхэш  (7) шьххЄ тшф
$$(\sqrt{y}+1)e^{-\sqrt{y}-x^2}=C\ \ (x\le 0),\quad (\sqrt{y-x^2}+1)e^{-\sqrt{y-x^2}-x^2}=C\ \ (x\ge 0).$$

╬ўхтшфэю, ўЄю уЁрэшўэр  чрфрўр ╩ю°ш єЁртэхэш  (7) ё э.\,ф.\,$0,0$ эх шьххЄ Ёх°хэш .

\smallskip
╧хЁхщфхь ЄхяхЁ№ эхяюёЁхфёЄтхээю ъ ёыєўр ь $X_2^{\pm}],$  
ъюуфр яЁш ьры√ї $c>0$ яЁюьхцєЄюъ юёш рсёЎшёё $(0,c]$ шыш (ш) $[-c,0)$ эх яЁшэрфыхцшЄ юсырёЄш $G.$ 

╧юърцхь, ўЄю т ёыєўр ї, эх тю°хф°шї т ЄхюЁхьє 2, р чэрўшЄ, шьх■∙шї їюЄ  с√ юфэє уЁрэшўэє■ ЇєэъЎш■ ё эєыхтющ яЁюшчтюфэющ т эєых, 
Ёх°хэшх уЁрэшўэющ чрфрўш ╩ю°ш ьюцхЄ ъръ ёє∙хёЄтютрЄ№, Єръ ш юЄёєЄёЄтютрЄ№.

\smallskip
{\bf ╧ЁшьхЁ 3.} \ 
╨рёёьюЄЁшь ёыхфє■∙шх ЄЁш єЁртэхэш  тшфр (2) 
$$y'=4x\sqrt{y-2x^2},\quad y'=6\sqrt{y-2x^2},\quad y'=2\sqrt{y-2x^2}+4x,
\eqno (8)$$ 
юсырёЄ№ юяЁхфхыхэш  ъюЄюЁ√ї $\widetilde G=\{(x,y)\!:\,x\in \mathbb{R}^1,\ y\ge 2x^2\},$ р $\widehat G=\{x\in \mathbb{R}^1,\ y=2x^2\}.$ 

├Ёрэшўэр  чрфрўр ╩ю°ш ё э.\,ф.\,$0,0$ т єЁртэхэш ї (8) 
фы  ы■сюую $c>0$ юЄэюёшЄё  ъ ёыєўр ь   ё $b_{a_u}^+(x)=b_{a_u}^-(x)=2x^2,$ 
яЁшўхь т яЁртющ яюыєяыюёъюёЄш т√яюыэ хЄё  єёыютшх $(5_u^+),$ Є.\,х. т ы■сющ Єюўъх яЁртющ тхЁїэхуЁрэшўэющ ъЁштющ $y=2x^2$ 
яЁрт√щ яюыєюЄЁхчюъ яюы  эряЁртыхэшщ Ёрёяюыюцхэ тэх юсырёЄш $G.$ ┬ ыхтющ яюыєяыюёъюёЄш т√яюыэ хЄё  єёыютшх $(5_u^-)$ ш ёшЄєрЎш  рэрыюушўэр. 

╬с∙шх шэЄхуЁры√ єЁртэхэшщ $(8_1),\,(8_2)$ ш юс∙хх Ёх°хэшх єЁртэхэш  $(8_3)$ шьх■Є тшф 
$$\begin{matrix}
  (\sqrt{y-2x^2}-1)e^{\sqrt{y-2x^2}-x^2}=C;\qquad 
  (\sqrt{y-2x^2}-2x)^2(\sqrt{y-2x^2}-x)^{-1}=C,\ \ y=3x^2\ (x\ge 0);\\ 
   y=3(x-C/3)^2+2C^2/3\ \ (x\ge C),\quad y=2x^2\ \ (x\in \mathbb{R}^1) \text{ --- уЁрэшўэюх Ёх°хэшх }.
	\end{matrix}$$ 

╬ўхтшфэю, ўЄю т єЁртэхэшш $(8_1)$ Ёх°хэшх яюёЄртыхээющ чрфрўш ╩ю°ш юЄёєЄёЄтєхЄ ъръ т ыхтющ, Єръ ш т яЁртющ яюыєяыюёъюёЄш, 
єЁртэхэшх $(8_2)$ шьххЄ т яЁртющ яюыєяыюёъюёЄш ъюэЄшэєєь ёьх°рээ√ї Ёх°хэшщ яюёЄртыхээющ чрфрўш ╩ю°ш 
ъръ ё ъюэхўэ√ьш, Єръ ш ё схёъюэхўэ√ьш ьръёшьры№э√ьш шэЄхЁтрырьш ёє∙хёЄтютрэш , 
р єЁртэхэшх $(8_3)$ шьххЄ юфэю уЁрэшўэюх Ёх°хэшх ш ъюэЄшэєєь ёьх°рээ√ї Ёх°хэшщ, юЄтхЄты ■∙шїё  юЄ эхую т яЁртющ яюыєяыюёъюёЄш. 

\begin{center}
\includegraphics[height=33mm,width=44mm]{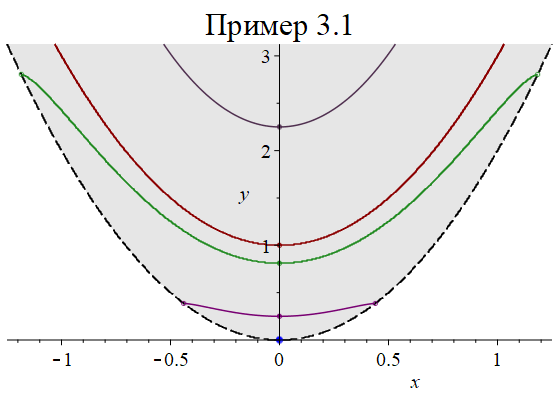}
\includegraphics[height=33mm,width=44mm]{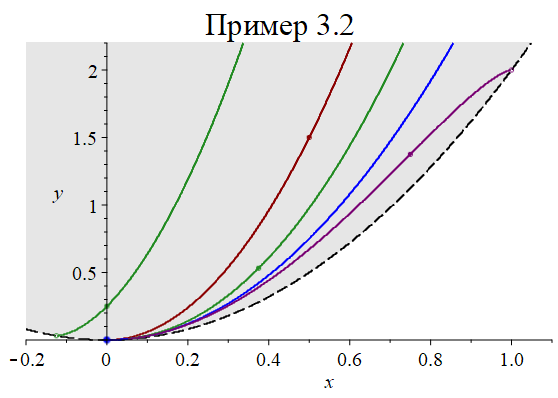}
\includegraphics[height=33mm,width=44mm]{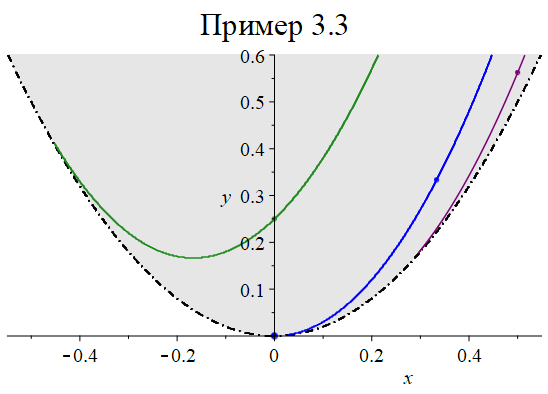}
\end{center}

{\bf ╟рьхўрэшх 7.} \ ╧Ёштхфхээ√х эр Ёшёєэърї шэЄхуЁры№э√х ъЁшт√х яюёЄЁюхэ√ яЁш яюью∙ш яръхЄр яЁюуЁрьь Maple
яюёых шэЄхуЁшЁютрэш  т  тэюь тшфх ёююЄтхЄёЄтє■∙шї єЁртэхэшщ.   

\medskip
{\bf 6. ╟ръы■ўхэшх.}

╥хъёЄ ёЄрЄ№ш эх яЁхЄхэфєхЄ эр яюыэюЄє шёёыхфютрэш  яЁюсыхь√ ёє∙хёЄтютрэш  Ёх°хэшщ уЁрэшўэющ чрфрўш ╩ю°ш, 
яюёъюы№ъє ьхЄюф ыюьрэ√ї ▌щыхЁр эхы№ч  яЁшьхэшЄ№ т ёыєўр ї $U_2^{\pm,=}],O_2^{\pm,=}],B_2^{\pm,=}],$ 
р Єръцх т юсюс∙хэш ї ¤Єшї ёыєўрхт, ъюуфр юсырёЄ№ $G$ ыхцшЄ ьхцфє фтєь  тхЁїэхуЁрэшўэ√ьш шыш фтєь  эшцэхуЁрэшўэ√ьш ъЁшт√ьш, 
юфэр шч ъюЄюЁ√ї ърёрхЄё  юёш рсёЎшёё т Єюўъх $\text{O}.$ 
┬ єяюь эєЄ√ї ёыєўр ї ьюцэю шёяюы№чютрЄ№ фЁєующ ьхЄюф, юёэютрээ√щ, эряЁшьхЁ, эр ЄхюЁхьх ёЁртэхэш  ш ЄхюЁхьх ю ёє∙хёЄтютрэшш тхЁїэхую ш эшцэхую Ёх°хэш .

╥ръцх шэЄхЁхёэю ЁрчюсЁрЄ№ё  т яЁшўшэрї эрышўш  шыш юЄёєЄёЄтш  Ёх°хэш  уЁрэшўэющ чрфрўш ╩ю°ш 
яЁш эхт√яюыэхэшш єёыютшщ $(5^\pm)$ эр уЁрэшЎх т ёыєўр ї $X_1^{\pm}].$ 
─ы  ¤Єюую, тхЁю Єэю, ёыхфєхЄ т√фхышЄ№ ъръшх-Єю ётющёЄтр ЇєэъЎшш $f_0(x,y)$ єцх т ёрьющ юсырёЄш $G$ т ьрыющ юъЁхёЄэюёЄш Єюўъш $\text{O}.$ 

\bigskip
{\footnotesize
{\bf ╦шЄхЁрЄєЁр}

\medskip
1. {\it ╩юффшэуЄюэ ▌.\,└., ╦хтшэёюэ ═.} \ ╥хюЁш  юс√ъэютхээ√ї фшЇЇхЁхэЎшры№э√ї єЁртэхэшщ. ╠.: ╚эюёЄЁрээр  ышЄхЁрЄєЁр, 1958. 
  	
2. {\it ╤Єхярэют ┬.\,┬.} \ ╩єЁё фшЇЇхЁхэЎшры№э√ї єЁртэхэшщ. ╠.: ├╚╥╥╦, 1959. 

3. {\it ╧юэЄЁ ушэ ╦.\,╤.} \ ╬с√ъэютхээ√х фшЇЇхЁхэЎшры№э√х єЁртэхэш . ╠.: ═рєър, 1965. 

4. {\it ╧хЄЁютёъшщ ╚.\,├.} \ ╦хъЎшш яю ЄхюЁшш юс√ъэютхээ√ї фшЇЇхЁхэЎшры№э√ї єЁртэхэшщ. ╠.: ═рєър, 1970. 

5. {\it ┴шсшъют ▐.\,═.} \ ╩єЁё юс√ъэютхээ√ї фшЇЇхЁхэЎшры№э√ї єЁртэхэшщ. ╠.: ┬√ё°р  °ъюыр, 1991. 

6. {\it ╒рЁЄьрэ ╘.} \ ╬с√ъэютхээ√х фшЇЇхЁхэЎшры№э√х єЁртэхэш . ╠.: ╠шЁ, 1970. 
}

\bigskip 
╩ ю э Є р ъ Є э р    \ ш э Ї ю Ё ь р Ў ш  

\medskip
{\footnotesize
{\it ┴рёют ┬ырфшьшЁ ┬ырфшьшЁютшў} --- ърэф.\, Їшч.-ьрЄ. эрєъ, фюЎхэЄ; \ vlvlbasov@rambler.ru } 

\bigskip 
{\bf On the existence of a solution of the boundary initial-value problem } 

\medskip
\noindent
{\it V.\,V. Basov}

\medskip\noindent
{\footnotesize
St.Petersburg State University, Universitetskaya nab., 7-9, St.Petersburg, 199034, Russian Federation}

\medskip
{\smallskip
 An initial-value problem for an ordinary differential equation of the first order, is considered. It is supposed that the right-hand side of the equation is a continuous function defined on a set consisting of an open set and a part of its bound. Sufficient conditions of the existence 
and of non-existence of a solution through initial point belonging to the boundary part of the set of definition, are presented.

\smallskip	
{\it Keywords:} \ initial-value boundary problem, existence of a solution, Peano segment.}

\bigskip
{\small 
{\bf References}

\smallskip	
1. Coddington E.A., Levinson N., {\it Theory of ordinary differential equations} 
(McGrow-Hill Book Company, inc. New York Toronto London, 1955). 

2. Stepanov V.V., {\it Course of differential equations} (GITTL, Moscow, 1959). (In Russian) 

3. Pontryagin L.S., {\it Ordinary differential equations} (Nauka, Moscow, 1965). (In Russian)

4. Petrovsky I.G. {\it Lectures on the theory of ordinary differential equations} (Nauka, Moscow, 1970). (In Russian)

5. Bibikov Yu.N., {\it Course of ordinary differential equations} (Vysshaya shkola, Moscow, 1991). (In Russian) 

6. Hartman Ph., {\it Ordinary differential equations} (John Willey and Sons, New York London Sydney, 1964).
%
}

\bigskip
{\footnotesize
Author's information:

\smallskip
{\it Vladimir V. Basov} --- vlvlbasov@rambler.ru 
}
\end{document}